 \documentclass[amsmath,amssymb,pra, showkeys]{revtex4-1}
 \def \LAA{false}

\usepackage{ifthen}
\usepackage{amssymb,amsmath}
\usepackage{graphicx}
\usepackage[format=hang, margin=1em]{subfig}
\usepackage{xcolor}

\def\ld{{\lambda}}
\def\dR{{\mathbb{R}}}
\def\dN{{\mathbb{N}}}

\def\dC{{\mathbb{C}}}

\def\conv{{\mbox{conv}}}

\def\diag{{\mbox{diag}}}
\def\trace{{\mbox{tr}}}

\def\fy{{\varphi}}

\newcommand{\beq}[2]{\begin{equation}\label{#1}#2}

\newlength{\figurewidth}
\ifthenelse{\equal{\LAA}{true}}{
 \setlength{\figurewidth}{2.2in}
}{
 \setlength{\figurewidth}{3.3in}
}

\begin{document}
\noindent

\ifthenelse{\equal{\LAA}{true}}
{
\begin{frontmatter}
\title{Numerical shadows: measures and densities\\
 on the numerical range
}
\author[virginia]{Charles F. Dunkl}
\author[iitis]{Piotr Gawron}
\author[guelph]{John A. Holbrook}
\author[iitis]{Zbigniew Pucha{\l}a}
\author[uj,cft]{Karol~{\.Z}yczkowski}
\address[virginia]{Department of Mathematics, University of Virginia, 
Charlottesville, VA 22904---4137, USA}
\address[iitis]{Institute of Theoretical and Applied Informatics, Polish Academy
of Sciences, Ba{\l}tycka 5, 44-100 Gliwice, Poland}
\address[guelph]{Department of Mathematics and Statistics, University of Guelph,
Guelph, Ontario, N1G 2W1, Canada}
\address[uj]{Instytut Fizyki im. Smoluchowskiego, Uniwersytet
Jagiello{\'n}ski, Reymonta 4, 30-059 Krak{\'o}w, Poland }
\address[cft]{Centrum Fizyki Teoretycznej, Polska Akademia Nauk, Aleja Lotnik{\'o}w 
32/44, 02-668 War\-sza\-wa, Poland}
%
}
{
%
\title{Numerical shadows: measures and densities on the numerical range}
\author{Charles F. Dunkl}
\affiliation{Department of Mathematics, University of Virginia, 
Charlottesville, VA 22904---4137, USA}
\author{Piotr Gawron}
\affiliation{Institute of Theoretical and Applied Informatics, Polish Academy
of Sciences, Ba{\l}tycka 5, 44-100 Gliwice, Poland}
\author{John A. Holbrook}
\affiliation{Department of Mathematics and Statistics, University of Guelph,
Guelph, Ontario, N1G 2W1, Canada}
\author{Zbigniew Pucha{\l}a}
\affiliation{Institute of Theoretical and Applied Informatics, Polish Academy
of Sciences, Ba{\l}tycka 5, 44-100 Gliwice, Poland }
\author{Karol \.Zyczkowski}
\affiliation {Instytut Fizyki im. Smoluchowskiego, Uniwersytet
Jagiello{\'n}ski,
ul. Reymonta 4, 30-059 Krak{\'o}w, Poland}
\affiliation{Centrum Fizyki Teoretycznej, Polska Akademia Nauk, Al.
Lotnik{\'o}w
32/44, 02-668 Warszawa, Poland}
%
}

\begin{abstract}

For any operator $M$ acting on an $N$-dimensional Hilbert space ${\cal H}_N$
we introduce its {\sl numerical shadow},  which is a probability measure 
on the complex plane supported by the numerical range of $M$.
The shadow of $M$ at point $z$ is defined as the probability that the 
inner product  $(Mu,u)$ is equal to $z$,
where $u$ stands for a random complex vector from ${\cal H}_N$,
satisfying $||u||=1$. In the case of $N=2$ the
numerical shadow of a non-normal operator 
can be interpreted as a shadow of a hollow sphere projected on a plane. 
A similar interpretation is provided also for higher dimensions.
For a hermitian $M$ its numerical shadow forms
a probability distribution on the real axis 
which is shown to be a one dimensional $B$-spline. 
In the case of a normal $M$ the numerical shadow 
corresponds to a shadow of a transparent solid simplex 
in ${\mathbb R}^{N-1}$ onto the complex plane. 
Numerical shadow is found explicitly for Jordan
matrices $J_N$, direct sums of matrices
and in all cases where the shadow is rotation invariant.
Results concerning the moments of shadow measures play an important role.
A general technique to study numerical shadow via the Cartesian decomposition
is described, and a link of the numerical shadow of an operator
to its higher-rank numerical range is emphasized.
\ifthenelse{\equal{\LAA}{true}}
{}{
{\\[.5cm]}
}
\end{abstract}

\ifthenelse{\equal{\LAA}{true}}
{
  \end{frontmatter}
}{
  \maketitle 
}


\section{Introduction}
The classical numerical range $W(M)$ of a
complex $N\times N$ matrix $M$ is the subset of $\dC$ defined by
$$W(M)=\{(Mu,u):u\in\dC^N, \|u\|=1\}.$$ This concept has a long
history and has proved a useful tool in operator theory and matrix
analysis as well as in more applied areas; for a nice account of
some of the lore of $W(M)$, see [GR1997]. Here we are concerned with
the measures or densities induced on $W(M)$ by various distributions
of the unit vector $u$; we use the term \textbf{numerical shadow} to
refer to such densities (see section 2 for our motivation in using
this terminology). Although it seems natural to study the numerical
shadow, this subject does not seem to have received much attention
until recently. The only earlier extended account (that we know of) is
the thesis [N1982]. This apparent neglect is not the only reason for
our attempts to understand the numerical shadow better; the numerical shadow also
plays an interesting role in quantum information theory.
Applications in this area are the focus of a companion paper
[DGHMP\.Z2], in preparation. Very recently a preprint by Gallay and Serre [GS2010]
has appeared, dealing also with mathematical aspects of the numerical shadow
(numerical measure, in their terminology). In several ways their development of the
subject is parallel to our own (as presented, for example, in [\.Z2009] and [H2010]).
In the present paper we stress those aspects of our work which are complementary to [GS2010],
particularly those based on the \textbf{moments} of the numerical shadow.
{\\[.5cm]}
In this paper we work mainly with
the ``uniform'' distribution of $u$ over the unit sphere $\Omega_N$
in $\dC^N\equiv\dR^{2N}$, ie the probability distribution on
$\Omega_N$ that is invariant under all orthogonal transformations of
$\dR^{2N}$. In some cases, however, the internal structure of $W(M)$
is better revealed through the use of other distributions, and these
are of special importance for the applications discussed in
[DGHMP\.Z2].
{\\[.5cm]}
Although the numerical shadow is in general
difficult to determine explicitly, this is possible in a number of
interesting cases. Methods based on identifying the moments of the
shadow measures are often effective. We also present, via the
figures, the results of numerical simulations that display shadow
densities in various other cases.
{\\[.5cm]}
In section 2 we treat the
simple situation occurring when $M$ is $2\times2$. Here we obtain a
``real--life'' shadow. In section 3 we discuss the analogous treatment
for $N\times N$ matrices $M$, ie we view the numerical shadow of $M$
as the image of an appropriate measure on the pure states $uu^*$
under a linear map $\fy_M$. In section 4 we show that in the case of
$N\times N$ normal $M$ the numerical shadow is the orthogonal
projection of a well-placed model of the ($N-1$)--dimensional
simplex (with the uniform density). Thus the density for the
numerical shadow is a 2--dimensional B-spline (1-dimensional in the
Hermitian case).
{\\[.5cm]}
 Section 5 studies the moments of
numerical shadows, yielding a key technique for the identification
and comparison of shadows. In section 6, criteria for the equality
of the numerical shadows of two matrices are obtained (in terms, for example, of
traces of words in the matrices and their adjoints). Evidently
equality occurs when the matrices are unitarily equivalent, but this
is not necessary (if $N>2$).
{\\[.5cm]}
In section 7, the numerical shadows
are found explicitly
for the Jordan nilpotents $J_N$.  Section 8 extends the techniques
developed in section 7 to obtain explicit densities for all rotation
invariant shadows.
{\\[.5cm]}
Section 9 introduces a useful view of the
numerical shadow in terms of the (Hermitian) components Re($M$) and
Im($M$) in the Cartesian decomposition of $M$, and the unitary
matrix linking those components. Several related aspects of the numerical shadow
are treated in that section, including a connection with the Radon transform.
Section 10 relates the numerical
shadow of a direct sum to the shadows of its summands.
Section 11 is concerned with numerical approximations of shadow densities
in terms of moments and Zernike expansions.
{\\[.5cm]}
Finally, in
section 12, we relate the numerical shadow of $M$ to the so--called
rank--$k$ numerical ranges $\Lambda_k(M)$. The theory and
applications of these ranges has been advanced vigorously since
their introduction only a few years ago as a tool in quantum information
theory (see for example [CK\.Z2006a, CK\.Z2006b, CHK\.Z2007, CGHK2008, W2008, LS2008, 
LPS2009, and GLW2010]).
One way to
describe $\Lambda_k(M)$ is that it consists of those points $(Mu,u)$
in $W(M)$ where $u$ may be chosen from the unit sphere in a whole
$k$--dimensional subspace of $\dC^N$. Thus it is natural to ask to
what extent $\Lambda_k(M)$ may be identified as a region of greater
density within the numerical shadow. Here the shadows corresponding
to \textbf{real} unit vectors $u$ play a role.
{\\[.5cm]}
Let us fix some notation. The algebra of complex $N\times N$ matrices is here denoted by
$M_N(\dC)$ or simply $M_N$. The adjoint or conjugate transpose of a matrix $M\in M_N$ is
denoted by $M^*$; we consider vectors $v$ in $\dC^N$ as column vectors and $v^*$ is the
conjugate transpose. Our inner product $(v,w)$ may be computed as $w^*v$. Recall that the unit sphere
in $\dC^N$ is denoted by $\Omega_N$, ie
\[
\Omega_N=\{u\in\dC^N:\|u\|=1\}.
\]
The uniform probability measure on $\Omega_N$ is denoted by $\mu$. Given $M\in M_N$, the notion of ``numerical shadow
of $M$'' is captured formally as the probability measure  $P_M$ on $W(M)$ such that
\[
P_M(S)=\mu\{u\in\Omega_N:(Mu,u)\in S\},
\]
for each Borel subset $S$ of $W(M)$. Equivalently, for any continuous function $g:W(M)\to\dC$ we have
\begin{equation}
{
\int_{W(M)} g(z)\,dP_M(z)=\int_{\Omega_N} g((Mu,u))\,d\mu(u).}
\label{E1.1}
\end{equation}

If $P_M$ has a probability density (with respect to planar measure in $\dC$) it is denoted by $f_M$. In those cases
where $f_M$ is rotation--invariant, ie $f_M(z)=f_M(|z|)$ we consider $f_M$ as a function of $r\in(0,w(M))$, where
$w(M)$ is the so--called numerical radius of $M$:
\[
w(M)=\max\{|z|:z\in W(M)\}.
\]
{\\[.5cm]}
Acknowledgements: Work by J. Holbrook was supported in part by an NSERC of Canada research grant.
Work by P. Gawron and Z.Pucha\l a was supported by the Polish Ministry of 
Science and Higher Education under the grant number N519 442339, 
while K.~{\.Z}yczkowski acknowledges support by the grant number 
  N202 090239.

\vskip 2.0cm  

\tableofcontents

\section{The $2\times2$ case; real--life shadows}

\begin{figure}[ht!]
\centering
\subfloat[][\label{fig:nonnormal}%
Numerical shadow of matrix nonnormal
$
\bigl(
\begin{smallmatrix} 
0 & 1 \\ 
0 & 0 
\end{smallmatrix}
\bigr)
$ resembles physical shadow cast by the hollow 
sphere made of transparent material when illuminated by a light source at infinity.] {%
\includegraphics[width=\figurewidth]{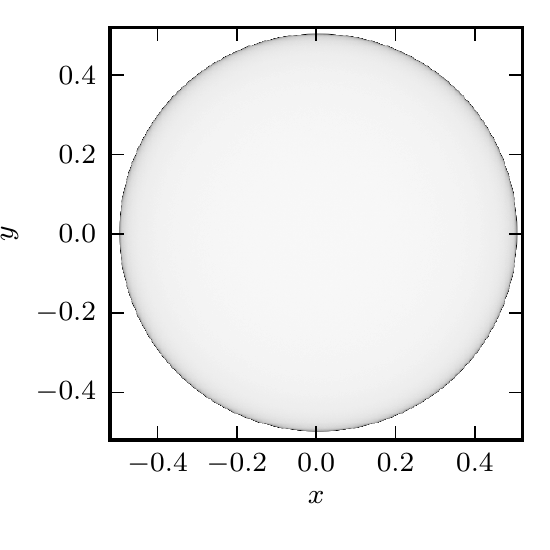}%
}%
\subfloat[][\label{fig:nonnormalcrossection}%
Cross-section of the shadow along the real axis.] {%
\includegraphics[width=\figurewidth]{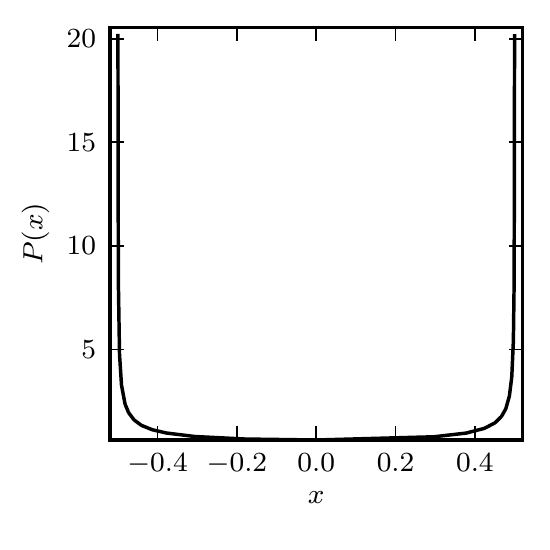}%
}%
\\
\subfloat[][\label{fig:nonnormal2}%
Numerical shadow of matrix nonnormal
$
\bigl(
\begin{smallmatrix} 
0 & 1 \\ 
0 & 1
\end{smallmatrix}
\bigr)
$ resembles physical shadow cast by the hollow 
sphere made of transparent material when illuminated by a light source at infinity, but with screen not perpendicular to the light rays.] {%
\includegraphics[width=\figurewidth]{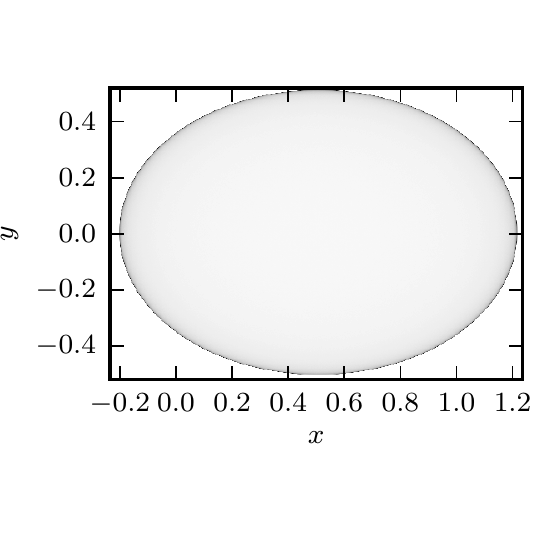}%
}%
\subfloat[][\label{fig:nonnormalcrossection2}%
Cross-section of the shadow along the real axis.] {%
\includegraphics[width=\figurewidth]{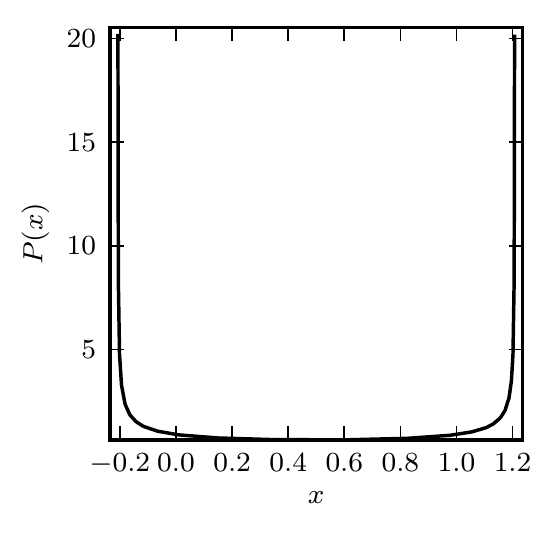}%
}%
\caption{ }
\end{figure}

Here we compute the shadow density for an arbitrary $2\times2$ matrix $M$.
Our method is to view this density as a real--life shadow of the Bloch sphere model for $\Omega_2$.
An equivalent result was obtained by Ng (see [N1982]) by a rather different method.
{\\[.5cm]}
The Bloch sphere model sees $W(M)$ as the genuine shadow of a hollow sphere made of infinitely thin semi--transparent uniform
material, where in general the light would fall obliquely on the (complex) plane. Of course, this situation cannot quite be
realized physically, but playing with a hollow plastic ball in bright sunlight may yield a good approximation.
It is well--known that when $M\in M_2$ the numerical range $W(M)$ is a filled ellipse with the eigenvalues of $M$ as
foci. The following proposition, visualized in Fig. \ref{fig:nonnormal2}, 
supplies further information in the form of an explicit shadow density.
{\\[.5cm]}
\textbf{Proposition 2.1:} Let $E$ be the filled ellipse formed by $W(M)$ and let $a$ and $b$ be the
lengths of the semimajor and semiminor axes of $E$;
then the shadow density is
\begin{equation}
\label{e1}{
\frac1{2\pi ab\sqrt{1-r^2}}}
\end{equation}
at every point on the elliptical curve bounding $rE$ ($0\leq r\leq1$).\\
\textbf{Proof:} Recall that
we assume that $u$ is chosen ``uniformly'' over
$\{u\in\dC^2:\|u\|=1\}$, ie according to the measure $\mu$ on $\Omega_2$. It is known that $|u_1|^2$ will then be
uniform in $[0,1]$. This a special case of the fact that $u$ uniform in $\Omega_N$ implies
$(|u_1|^2,|u_2|^2,\dots,|u_N|^2)$ has the uniform distribution in the $N-1$--dimensional simplex, see  \cite{Z1999,BZ2006}.
%
It is also important to note that for fixed $u_1$ the relative phase of $u_2$ is $e^{i\theta}$ where $\theta$ is uniform in
$[0,2\pi]$. We have
\[
uu^*=\begin{bmatrix}|u_1|^2&u_1\overline{u_2}\\ \overline{u_1}u_2&|u_2|^2\end{bmatrix} =
\]\[
\frac12\begin{bmatrix}1&0\\0&1\end{bmatrix}+\frac12\begin{bmatrix}2|u_1|^2-1&2u_1\overline{u_2}\\
2\overline{u_1}u_2&1-2|u_1|^2\end{bmatrix}.
\]
Let $z=2|u_1|^2-1$; then $z$ is uniform in $[-1,1]$ and
$2u_1\overline{u_2}=\sqrt{1-z^2}e^{-i\theta}=x-iy$, so that (recalling Archimedes) $(x,y,z)$ is uniform on the unit sphere. This is
one way to see that the corresponding distribution on the Bloch sphere $\{uu^*:u\in\Omega_2\}$ is uniform.
{\\[.5cm]}
Following Davis [D1971], we compute $(Mu,u)$ as
\[
\trace(u^*Mu)=\trace(Muu^*)=
\trace(\frac12 M + \frac12 M\begin{bmatrix}z&x-iy\\x+iy&-z\end{bmatrix}).
\]
For convenience we take $M=\begin{bmatrix}1&q\\0&-1\end{bmatrix}$ with $q\geq0$. This is a harmless normalization, achieved
via translation and rotation of $M$ (which respects the distribution of $(Mu,u)$) and unitary similarity (which leaves
the distribution unchanged). Then $(Mu,u)=z+(q/2)(x+iy)=(z+bx,by)$ in the complex plane, with $b=q/2$. Consider the
region $E$ bounded by the ellipse centred at $(0,0)$ with horizontal semiaxis of length $a=\sqrt{1+b^2}$ and
vertical semiaxis of length $b$. Given $r\in[0,1]$, $(Mu,u)$ lies in $rE$ iff
\[
\frac{(z+bx)^2}{1+b^2}+\frac{(by)^2}{b^2}\leq r^2.
\]
A calculation verifies that this is equivalent to $((x,y,z)\cdot(1/a,0,-b/a))^2\geq 1-r^2$, saying that $(x,y,z)$ lies
on either of the spherical caps of the unit sphere that are symmetrical about the axis determined by $(1/a,0,-b/a)$ and have radius $r$.
According to Archimedes (or the related formulas found in calculus texts) the relative area of these caps is
$1-\sqrt{1-r^2}$. Hence the probability
\[
P((Mu,u)\in rE)=1-\sqrt{1-r^2}.
\]
To find the corresponding planar density, we first observe that the region $(r+\Delta r)E\setminus rE$ corresponds to symmetrical
rings bordering the spherical caps mentioned above. Thus the planar density will be constant on the ellipse bounding $rE$. Its
value there is then given by
\[
\lim_{\Delta r\to0}\frac{(1-\sqrt{1-(r+\Delta r)^2})-(1-\sqrt{1-r^2})}{\pi ab(r+\Delta r)^2 - \pi abr^2}=\frac1{2\pi ab\sqrt{1-r^2}}.
\]
QED
{\\[.5cm]}
This result is equivalent to the formula for the density obtained by Ng (see [N1982], page 67). He computes the density
of $(Mu,u)$ at $(x,y)$ in the ellipse as
\[
   p(x,y)= \frac1{2\pi ab\sqrt{1-(x^2/a^2+y^2/b^2)}}.
\]
His method seems unrelated to the Bloch sphere approach worked out above.

\section{Numerical shadows as linear images of the pure quantum states}
It is clear that the argument of section 2 can be extended in part to cases where $N>2$.
The Bloch sphere is replaced by the set of density matrices representing pure quantum states:
\[
PQS_N=\{uu^*:u\in\Omega_N\}.
\]
Just as before, for any $u\in\Omega_n$ and $M\in M_N$ we have
\[
(Mu,u)=\trace(u^*Mu)=\trace(Muu^*),
\]
so that $W(M)$ is the image of $PQS_N$ under the linear map $\fy_M:M_N\to\dC$ defined by
\[
\fy_M(X)=\trace(MX).
\]
Since each $X\in PQS_N$ is Hermitian we may also write
\[
\fy_M(X)=\trace(MX^*)=(M,X)_F,
\]
where $(\cdot,\cdot)_F$ is the Frobenius inner product on $M_N$.
{\\[.5cm]}
Thus the numerical shadow of $M$ may be viewed as the measure on $W(M)$ induced by applying the linear
map $\fy_M$ to the fixed measure $\nu$ on $PQS_N$ that corresponds to the uniform $\mu$ on $\Omega_N$.
As $M$ varies the resulting numerical shadows may be regarded as a tomographic study of the measure $\nu$.
Thus such detailed information as we have about numerical shadows (see section 9, for example, or [GS2010])
reveals much about the structure of $\nu$  on $PQS_N$.

\section{The Hermitian and normal cases: B--splines}

\begin{figure}[ht!]
\centering
\subfloat[][\label{fig:hermitian3}%
Shadow of matrix $\diag{(0, 1, 3)}$ is a spline function of degree one.] {%
\includegraphics[width=\figurewidth]{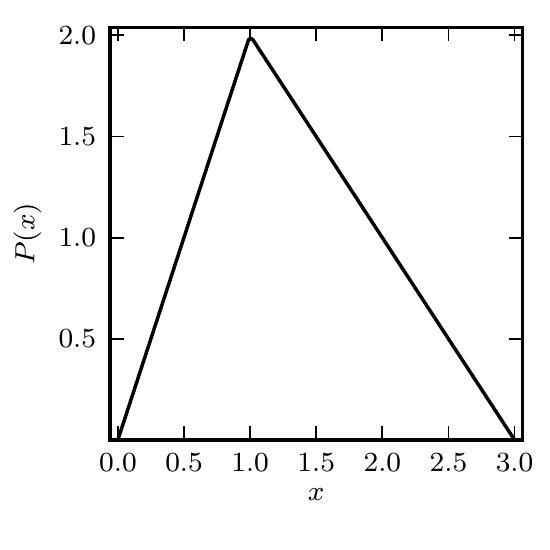}%
}%
\subfloat[][\label{fig:hermitian4}%
Shadow of matrix $\diag{(0, 1, 3, 5)}$ is a spline function of degree two.] {%
\includegraphics[width=\figurewidth]{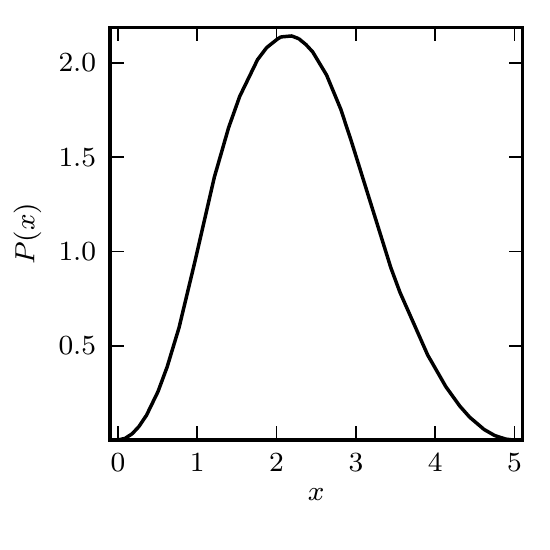}%
}%
\caption{Probability density function of hermitian matrices of dimensions 3 and 4.}
\label{fig:herm}
\end{figure}

\begin{figure}
\centering
\subfloat[][\label{fig:unitary3}%
Shadow of the matrix $\diag{(e^{2/3 i \pi },e^{-2/3 i \pi},1)}$ forms an 
uniform distribution supported by the triangle spanning its eigenvalues] {%
\includegraphics[width=\figurewidth]{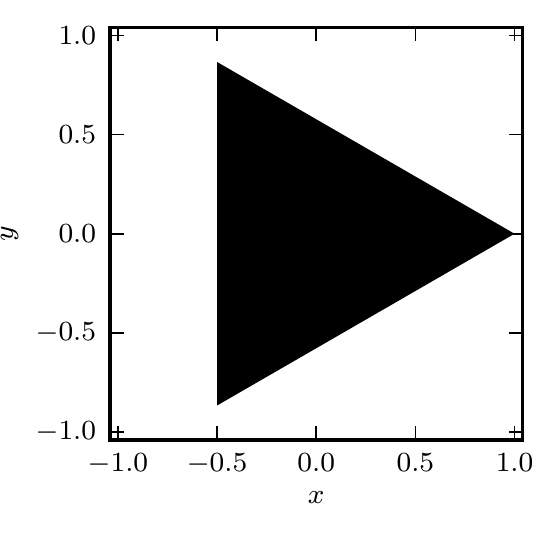}%
}%
\subfloat[][\label{fig:unitary4}%
Shadow of the matrix $\diag{(1, i, -1, -i)}$ forms a regular pyramid.] {%
\includegraphics[width=\figurewidth]{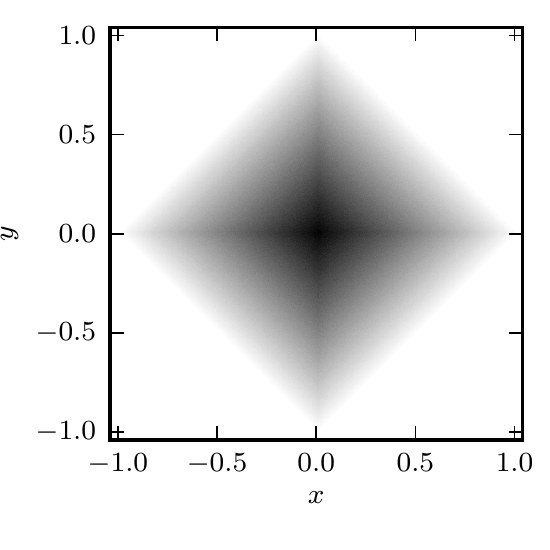}%
}%

\subfloat[][\label{fig:unitary5}%
Shadow of the matrix $\diag{(e^{2/5 i \pi }, e^{4/5 i \pi }, e^{6/5 i \pi }, e^{8/5 i \pi }, 1)}$.] {%
\includegraphics[width=\figurewidth]{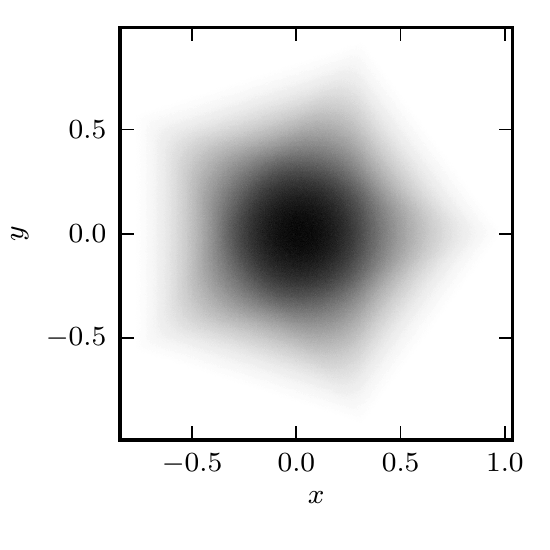}%
}%
\caption{Probability density function of unitary matrices of dimensions 3, 4 and 5.}
\label{fig:norm}
\end{figure}

The standard $N$--simplex $\Delta_N$ is defined by
\[
\Delta_N=\{r\in\dR^N: \mbox{all }r_k\geq0\mbox{ and }\sum_{k=1}^N r_k=1\}.
\]
This is an ($N-1$)--dimensional convex subset of $\dR^N$. We say $r$ is ``uniformly distributed'' over $\Delta_N$
to mean that $r$ is uniform with respect to normalized ($N-1$)--dimensional Lebesgue measure $vol_{N-1}$ on $\Delta_N$.
{\\[.5cm]}
\textbf{Lemma 4.1:} Let $P_A$ be the shadow measure of  a normal matrix $A\in M_N$. Then for any Borel subset $B$
of the plane
\[
P_A(B)=\mbox{Prob}\{r^*\ld\in B\},
\]
where $\ld=(\ld_1,\dots,\ld_N)^*$ is the spectrum of $A$ and $r$ is uniformly distributed over the standard
$N$--simplex $\Delta_N$.\\
\textbf{Proof:} Since $A$ is normal, it is unitarily similar to $\diag(\ld)$. As $\mu$ is invariant under unitary transformations we have
\[
P_A(B)=\mbox{Prob}\{(\diag(\ld)u,u)\in B:u\in\Omega_N\}=\mbox{Prob}\{(\sum_{k=1}^N\ld_k|u_k|^2\in B:u\in\Omega_N\}.
\]
It is known  (see e.g. \cite{Z1999,BZ2006})
that if $u$ is uniform over $\Omega_N$ (ie distributed according to $\mu$)
then
\[
r=(|u_1|^2,\dots,|u_N|^2)
\]
is uniform over the simplex $\Delta_N$. QED
{\\[.5cm]}
With $\ld$ as above, write $\ld=x+iy$ where $x,y\in\dR^N$. Assume first that $x$ and $y$ are independent vectors.
Let $W$ be a real matrix with $x,y$ as its first two columns and with columns 3 to $N$ forming an orthonormal
basis for $\{x,y\}^\perp$. For any vector $v\in\dR^N$,
\[
v^*\ld=a+ib\,\,\,(a,b\in\dR)\quad\Leftrightarrow\quad(W^*v)_1=a\mbox{ and }(W^*v)_2=b.
\]
Thus for Borel $B\subset\dC\equiv\dR^2$
\[
\mbox{Prob}\{r^*\ld\in B\}=vol_{N-1}\{r\in\Delta_N:((W^*r)_1,(W^*r)_2)\in B\},
\]
and, in view of Lemma 4.1, the density for $P_A$ at $(a,b)\in W(A)$ is
\[
vol_{N-3}\{r\in\Delta_N:(W^*r)_1=a,(W^*r)_2=b\}=vol_{N-3}\{v\in W^*(\Delta_N):v_1=a,v_2=b\}.
\]
Let us recall the definition of an $s$--dimensional B--spline (from [dB1976]).
{\\[.5cm]}
\textbf{Definition 4.2:} Let $\sigma$ be a nontrivial simplex in $\dR^{s+k}$. On $\dR^s$ we
define the B--spline of order $k$ from $\sigma$ by
\[
M_{k,\sigma}(x_1,\dots,x_s)=vol(\sigma\cap\{v\in\dR^{s+k}:v_j=x_j\,\,(j=1,2,\dots,s)\}).
\]
Using this terminology we may summarize our results as follows.
{\\[.5cm]}
\textbf{Proposition 4.3:} The numerical shadow of an $N\times N$ normal matrix with eigenvalues
$\ld\in\dC^N$ having linearly independent real and imaginary parts has as density a 2--dimensional
B--spline $M_{N-2,\sigma}(a,b)$, where the simplex $\sigma=W^*(\Delta_N)$ with some $W$ chosen as above.
{\\[.5cm]}
\textbf{Remark 4.4:} In the case where the real and imaginary parts of $\ld$ are dependent (as when $A$ is
Hermitian), it is easy to see that the numerical shadow is 1--dimensional (a line segment, in fact)
with density given by a 1--dimensional B--spline (compare [dB1976, Lemma 9.1]).
As examples, the shadows of Hermitian matrices of size $N=3$ and $N=4$
are shown in Fig. \ref{fig:herm}.
{\\[.5cm]}
\textbf{Remark 4.5:} This observation for the Hermitian case was worked out in detail by Ng in [N1982].
He also made the right conjecture regarding the normal case. In a sense, the normal case was earlier understood by statisticians
studying the distribution of quadratic forms; see for example [A1971, chapter 6]; Anderson points out that
some of the relevant ideas go back to von Neumann in the 40's. Anderson seems to discuss only real quadratic forms; thus the
normal case corresponds to roots of multiplicity 2.
{\\[.5cm]}
\textbf{Remark 4.6:} In view of Proposition 4.3, the theory of B--splines may be applied to see that normal
shadow densities are piecewise polynomial functions of two variables
in the normal case - see Fig. \ref{fig:norm}
and of a single variable in the Hermitian case.
The latter case is analyzed in some detail in Sec. 9.
A thorough analysis of the B--spline shadow 
densities for normal matrices is also provided in [GS2010].

\section{Moments of the numerical shadow}
We denote the moments of the numerical shadow of $A\in M_N$ by
\[
\nu_{jk}(A)=\int_{W(A)} z^j\overline{z}^k\,dP_A(z).
\]
Note that, since the polynomials in $z$ and $\overline{z}$ are uniformly dense in the continuous functions on
$W(A)$, these moments determine the numerical shadow uniquely. Moreover, in view of (\ref{E1.1}), we have
\beq{E6.1}{
\nu_{jk}(A)
=\int_{\Omega_N} (Au,u)^j(\overline{(Au,u)})^k\,d\mu(u)
=\int_{\Omega_N} (Au,u)^j(A^*u,u)^k\,d\mu(u).}
\end{equation}
Given $\ld\in\dC^N$ and a multi--index $\alpha\in\dN_0^N$ (where $\dN_0=\{0,1,2,\dots\}$), we use the following notation:
$\ld^\alpha=\ld_1^{\alpha_1}\ld_2^{\alpha_2}\dots\ld_N^{\alpha_N}$, $|\alpha|=\sum_1^N \alpha_k$,
$\alpha!=\alpha_1!\alpha_2!,\dots\alpha_N!$. We also use the Pochhammer symbol or shifted factorial
$(x)_n=\prod_{j=1}^n(x+j-1)$; by convention $(x)_0=1$.
{\\[.5cm]}
The effective evaluation of the moments $\nu_{jk}(A)$ depends on the following proposition.
{\\[.5cm]}
\textbf{Proposition 5.1:} Given $A\in M_N$, let $\ld\in\dC^N$ list the eigenvalues of $A$ repeated according to multiplicity.
Then
\beq{E6.2}{
\int_{\Omega_N} (Au,u)^n\,d\mu(u)=\frac{n!}{(N)_n} h_n(\ld),}
\end{equation}
where $h_n(\ld)$ is the complete symmetric polynomial of degree $n$, ie
\[
h_n(\ld)=\sum_{\alpha\in\dN_0^N,\,\,|\alpha|=n} \ld^\alpha.
\]
\textbf{Proof:} Given multi--indices $\alpha,\beta$, let
\beq{E6.3}{
Q(\alpha,\beta)=\int_{\Omega_N} u^\alpha (\overline{u})^\beta\,d\mu(u),}
\end{equation}
where the conjugation $\overline{u}$ is applied entrywise. Since $\mu$ is invariant under the unitary map
$u\to(e^{i\theta}u_1,u_2,\dots,u_N)^t$, we have $Q(\alpha,\beta)=e^{i(\alpha_1-\beta_1)\theta}Q(\alpha,\beta)$ for each
real $\theta$. Hence $Q(\alpha,\beta)=0$ unless $\alpha_1=\beta_1$. Similarly for the other components, so that
$Q(\alpha,\beta)=0$ unless $\alpha=\beta$. More work is required to evaluate $Q(\alpha,\alpha)$:
\beq{E6.4}{
Q(\alpha,\alpha)=\int_{\Omega_N}|u|^{2\alpha}\,d\mu(u)=\frac{\alpha!}{(N)_{|\alpha|}},}
\end{equation}
where $|u|=(|u_1|,|u_2|,\dots,|u_N|)^t$. A convenient trick here is to consider
\[
I=\int_{\dR^{2N}}e^{-\sum_1^N(x_k^2+y_k^2)}\prod_1^N(x_k^2+y_k^2)^{\alpha_k}\,\,dx\,dy.
\]
As a product of Gamma--integrals we obtain $I=\pi^N\alpha!$. Integrating first over $r\Omega_N$, with $r^2=\sum_1^N(x_k^2+y_k^2)$,
then over $0<r<\infty$, we find that $I=\frac12|S^{2N-1}|(N+|\alpha|-1)!Q(\alpha,\alpha)$, where $|S^{2N-1}|$ denotes the $(2N-1)$--dimensional
area of $\Omega_N$. Since $Q(\vec0,\vec0)=1$, the formula (\ref{E6.4}) follows.
{\\[.5cm]}
We may assume $A$ is in the Schur upper--triangular form, since this is obtained via a unitary similarity and $\mu$ is
invariant under unitary transformations on $\dC^N$. Thus $A_{jj}=\ld_j$ (some listing of the eigenvalues of $A$, with
multiplicity) and
\[
(Au,u)=\sum_{j=1}^N\ld_j|u_j|^2+\sum_{j>i}A_{ij}u_j\overline{u_i}.
\]
Aside from $(\sum_j \ld_j|u_j|^2)^n$, the terms of $(Au,u)^n$ are scalar multiples of expressions of the form
\[
\prod_{k=1}^a |u_{\ell_k}|^2\prod_{k=1}^b u_{j_k}\overline{u_{i_k}},
\]
where $b\geq1$ and each $j_k>i_k$. Such an expression has the form $u^\gamma(\overline{u})^\gamma u^\alpha(\overline{u})^\beta$
where, using $e(j)$ as a temporary notation for the multi--index with 1 in the $j$--th position and 0's elsewhere,
\[
\alpha=\sum_{k=1}^b e(j_k),\quad\beta=\sum_{k=1}^b e(i_k).
\]
Clearly $\alpha_1=0$ and for some first $k_0>1$ we have $\alpha_{k_0}>0$; hence $\beta_k>0$ for some $k<k_0$, so that
$\alpha\neq\beta$. Thus $Q(\gamma+\alpha,\gamma+\beta)=0$ so that such terms make no contribution to the integral over
$\Omega_N$. It follows that
\[
\int_{\Omega_N} (Au,u)^n\,d\mu(u)=\int_{\Omega_N}(\sum_{j=1}^N\ld_j|u_j|^2)^n\,d\mu(u).
\]
Using the multinomial formula and (\ref{E6.4}), this integral is
\[
\sum_{\alpha\in\dN_0^N,\,\,|\alpha|=n}\frac{n!}{\alpha!} \ld^\alpha Q(\alpha,\alpha)=
\sum_{\alpha\in\dN_0^N,\,\,|\alpha|=n}\frac{n!}{\alpha!} \ld^\alpha \frac{\alpha!}{(N)_n}
\]
\[
=\frac{n!}{(N)_n}\sum_{\alpha\in\dN_0^N,\,\,|\alpha|=n}\ld^\alpha=\frac{n!}{(N)_n} h_n(\ld).
\]
QED
{\\[.5cm]}
It will be convenient to use the notation $\ld(A)$ to denote any listing of the eigenvalues of $A$, repeated according
to multiplicity. Applying (\ref{E6.2}) with $A$ replaced by $sA+tA^*$ ($t,s$ real) and recalling (\ref{E6.1}) we obtain
\begin{equation}
\label{E6.5}
\sum_{j=0}^n\binom{n}j 
s^j t^{n-j}\nu_{j,n-j}(A)=\frac{n!}{(N)_n} h_n(\ld(sA+tA^*)).
\end{equation}
Moreover, the RHS of (\ref{E6.5}) may be evaluated in terms of traces of words in $A$ and $A^*$, using known relations \cite{Mac1995}
among $h_n(\ld)$, the power sums 
\[
p_j(\ld)=\sum_{k=1}^N \ld_k^j
\]
(these equal $\trace(A^j)$ if $\ld=\ld(A)$), and the elementary symmetric polynomials
\[
e_j(\ld)=\sum_{1\leq i_1<i_2<\dots<i_j\leq N} \ld_{i_1}\ld_{i_2}\dots\ld_{i_j}
\]
(note that $e_j(\ld)=0$ if $j>N$; by convention $e_0(\ld)=1$). For $n\geq1$ we have
\beq{E6.6}{
e_0h_n-e_1h_{n-1}+e_2h_{n-2}\dots\pm e_nh_0=0}
\end{equation}
(by convention $h_0(\ld)=1$) and
\beq{E6.7}{
ne_n=p_1e_{n-1}-p_2e_{n-2}+\dots\pm p_ne_0.}
\end{equation}

For example, $h_1=p_1$ so that (\ref{E6.5}) implies
\begin{eqnarray*}
t\nu_{0,1}(A)+s\nu_{1,0}(A) &=& \frac1N p_1(\ld(sA+tA^*)) = \frac1N\trace(sA+tA^*)\\
&=&\frac1N(t\,\,\trace(A^*)+s\,\,\trace(A)).
\end{eqnarray*}

Hence
\[
\nu_{1,0}(A)=\frac1N\trace(A)\mbox{ and }\nu_{0,1}(A)=\frac1N\trace(A^*).
\]
Likewise $1h_2=e_1h_1-e_2h_0=e_1h_1-\frac12(p_1e_1-p_2e_0)=p_1^2-\frac12p_1^2+\frac12p_2=\frac12(p_2+p_1^2)$,
so that (\ref{E6.5}) implies

\begin{eqnarray*}
 && t^2\nu_{0,2}(A)+2ts\nu_{1,1}(A)+s^2\nu_{2,0}(A) = \frac2{N(N+1)}(\frac12\trace(sA+tA^*)^2+\frac12\trace^2(sA+tA^*)) \\
 && =\frac1{N(N+1)}\big(t^2(\trace(A^*)^2+\trace^2A^*)+2ts(\trace(AA^*)+\trace(A)\trace(A^*)) \\
 && \phantom{ = } +s^2(\trace(A^2)+\trace^2A)\big) 
\end{eqnarray*}
(using $\trace(A^*A)=\trace(AA^*)$). Thus we have
\[
\nu_{2,0}(A)=\frac{\trace(A^2)+\trace^2A}{N(N+1)},\,\,\nu_{1,1}(A)=\frac{\trace(AA^*)+\trace(A)\trace(A^*)}{N(N+1)},
\]
and
\[
\nu_{0,2}(A)=\frac{\trace(A^*)^2+\trace^2A^*}{N(N+1)}.
\]
{\\[.5cm]}
Similarly we find that $h_3=\frac13p_3+\frac12p_1p_2+\frac16p_1^3$, so that (\ref{E6.5}) implies
\begin{eqnarray*}
&&  t^3\nu_{0,3}(A)+3t^2s\nu_{1,2}(A)+3ts^2\nu_{2,1}(A)+s^3\nu_{3,0}(A) \\ 
&& =\frac6{N(N+1)(N+2)}\big(\frac13\trace(sA+tA^*)^3+\frac12\trace(sA+tA^*)\trace(sA+tA^*)^2\\
&& \phantom{ = }+\frac16\trace^3(sA+tA^*)\big)\\
&& =\frac1{N(N+1)(N+2)}\big(t^3(2\trace(A^*)^3+3\trace(A^*)^2\trace(A^*)+\trace^3A^*)\\
&& \phantom{ = } + t^2s(6\trace(A(A^*)^2)+6\trace(A^*)\trace(AA^*)\\
&& \phantom{ = } + 3\trace(A)\trace(A^*)^2+3\trace(A)\trace^2A^*)+ts^2\dots\big),
\end{eqnarray*}

where again the cyclicity of the trace plays a role.
{\\[.5cm]}
From these calculations we obtain
\begin{eqnarray*}
\nu_{3,0}(A)&=&\frac1{N(N+1)(N+2)}\bigl(
   2\trace(A^3)+3\trace(A^2)\trace(A)+\trace^3A\bigr), \\
\nu_{2,1}(A)&=&\frac1{N(N+1)(N+2)}\bigl( 2\trace(A^2A^*)+2\trace(A)\trace(AA^*) \\
& & +\trace(A^2)\trace(A^*)+\trace^2A\,\,\trace(A^*)\bigr),
\end{eqnarray*}
and $\nu_{1,2}(A),\nu_{0,3}(A)$ by interchanging the roles of $A$ and $A^*$.
{\\[.5cm]}
What is not clear from the approach above is the important fact that all the moments $\nu_{j,k}(A)$ are polynomials
in the traces of $A,A^*$ words of length at most $N$. One way to see this is to note that (\ref{E6.6}) and (\ref{E6.7})
allow us to express $h_n$ for $n>N$ in terms of $p_1,p_2,\dots,p_N$. For example, when $N=3$ we find that
$h_4=\frac23p_1p_3+\frac1{12}p_1^4+\frac14p_2^2$.
{\\[.5cm]}
In general then we need only compute $\trace(sA+tA^*)^k$ for $k\leq N$, and therefore (in view of the noncommutative binomial
formula) we need only compute $\trace(P_{k,j}(A,A^*))$ for $k\leq N$, where
\[
P_{k,j}(x,y)=\sum_{S\subseteq\{1,2,\dots,k\},\,\#(S)=j} \{z_1z_2\dots z_k
\mbox{ where }z_i=x \mbox{ if }i\in S,\,z_i=y \mbox{ if }i\not\in S\}.
\]
{\\[.5cm]}
We summarize this discussion in the following proposition.
{\\[.5cm]}
\textbf{Proposition 5.2:} Given $A\in M_N$, all the moments $\nu_{i,j}(A)$ of the shadow measure (and therefore
the measure $P_A$ itself) are determined by the traces of $(sA+tA^*)^k$ (as polynomials in $s,t$) for $k\leq N$.
Thus they are determined by the values of
\[
\trace\bigl( P_{k,j}(A,A^*)\bigr)
\]
for $k\leq N$.
{\\[.5cm]}
\textbf{Remark:} The cyclicity of the trace (ie $\trace(AB)=\trace(BA)$) reduces $\trace\bigl( P_{k,j}(A,A^*)\bigr)$ to a single
term when $k\leq3$ but this is not always the case. For example
\[
\trace\bigl( P_{4,2}(A,A^*)\bigr)=4\trace\bigl( A^2(A^*)^2\bigr)+2\trace(AA^*)^2.
\]
{\\[.5cm]}
The information about moments $\nu_{j,k}(A)$ that is provided by (\ref{E6.5}) may also be encoded in the series
\begin{eqnarray*}
S(A,s,t,q) &=& \sum_{n=0}^\infty q^n\frac{(N)_n}{n!}(\sum_{j=0}^n\binom{n}{j}s^jt^{n-j}\nu_{j,n-j}(A))\\
&=& \sum_{n=0}^\infty q^n h_n(\ld(sA+tA^*))
\end{eqnarray*}
(absolutely convergent for small real $s,t,q$).
{\\[.5cm]}
The following proposition provides powerful alternative forms for this series.
{\\[.5cm]}
\textbf{Proposition 5.3:} Given $A\in M_N$ we have (for sufficiently small $s,t,q$)
\beq{E6.8}{
S(A,s,t,q)=\mbox{det}^{-1}\bigl(I-q(sA+tA^*)\bigr)}
\end{equation}
and
\beq{E6.9}{
S(A,s,t,q)=\Bigl(\sum_{k=0}^N(-q)^k e_k\bigl(\ld(sA+tA^*)\bigr)\Bigr)^{-1}.}
\end{equation}
\textbf{Proof:} These follow from the identities
\begin{eqnarray*}
\sum_{n=0}^\infty q^nh_n(\ld) &=& \prod_{j=1}^N(1+q\ld_j+q^2\ld_j^2+\dots) \\ 
&=&\bigl(\prod_{j=1}^N(1-q\ld_j)\bigr)^{-1}=\bigl(\sum_{k=1}^N(-q)^ke_k(\ld)\bigr)^{-1}.
\end{eqnarray*}
QED
{\\[.5cm]}
\textbf{Remark 5.4:} In view of (\ref{E6.8}), the shadow measure $P_A$ is completely determined by
det$\bigl(I-(sA+tA^*)\bigr)$ as a polynomial in $s,t$ ($q$ may be absorbed into $t,s$).
{\\[.5cm]}
\textbf{Remark 5.5:} In view of (\ref{E6.7}), the relation (\ref{E6.9}) provides another viewpoint
on Proposition 5.2.

\section{Criteria for equality of numerical shadows}

Given $A,B\in M_N$, we have seen in the last section that $P_A=P_B$\\
iff
\beq{E7.1}{
\trace(sA+tA^*)^k=\trace(sB+tB^*)^k}
\end{equation}
(as polynomials in $s,t$) for all $k\leq N$\\
iff
\beq{E7.2}{
\trace\bigl(P_{k,j}(A,A^*)\bigr)=\trace\bigl(P_{k,j}(B,B^*)\bigr)}
\end{equation}
for all $j\leq k\leq N$\\
iff
\beq{E7.3}{
\det(I-(sA+tA^*))=\det\bigl(I-(sB+tB^*)\bigr)}
\end{equation}
(as polynomials in $s,t$).
{\\[.5cm]}
Since the uniform measure $\mu$ on $\Omega_N$ is invariant under unitary transformations, $P_A=P_{U^*AU}$
for any unitary $U$. It is natural, therefore, to ask whether the numerical shadow $P_A$ determines $A$ up to
unitary similarity. This is the case for $A\in M_2$, for example, since the ellipse $E=W(A)$, just as a set, determines
an upper--triangular form for $A$: $A$ is unitarily similar to $\begin{bmatrix}\alpha&s\\0&\beta\end{bmatrix}$ where
the eigenvalues $\alpha,\beta$ are the foci of $E$ and $s$ is the length of the minor axis of $E$. The answer is ``yes" also
for normal matrices $A$ since the eigenvalues are determined by $P_A$ in that case (see section 4).
{\\[.5cm]}
More generally, however, the answer is ``no", on several levels. First of all, the measure $\mu$ is also invariant under
any orthogonal transformation of $\dR^{2N}\equiv\dC^N$, so that, in particular, $d\mu(u)=d\mu(\overline{u})$. Thus $A$ and
its transpose $A^t$, though they are not usually unitarily similar, always have the same numerical shadow:
\[
(Au,u)=u^*Au=(u^*Au)^t=u^tA^t(u^*)^t=(\overline{u})^*A^t\overline{u}=(A^t\overline{u},\overline{u}).
\]
In fact, the maps $A\mapsto U*AU$ and $A\mapsto U^*A^tU$, are the only linear maps on $M_N$ that preserve the
numerical shadow, since they are the only linear maps that preserve the numerical range as a set (see C.--K. Li's survey
[L2001]).
{\\[.5cm]}
Moreover, particular pairs $A,B$ may have the same numerical shadow without being related by unitarily similarity or transpose.
This phenomenon is somewhat clarified by comparing the trace criterion (\ref{E7.2}) for $P_A=P_B$ with the analogous
criteria for unitary similarity. In a 1940 paper [S1940] Specht observed that $A$ and $B$ are unitarily similar iff
\[
\trace\bigl(w(A,A^*)\bigr)=\trace\bigl(w(B,B^*)\bigr)
\]
for all two--variable words $w(\cdot,\cdot)$. Since then much work has been done with the aim of limiting the set
of words required  in Specht's criterion when matrices of a given size are involved. In [DJ2007] Djokovi\'c and
Johnson provide a welcome
account of recent results in this direction. In particular, the following result (see Theorem 2.4 in {DJ2007})
may be compared with (\ref{E7.2}).
{\\[.5cm]}
\textbf{Proposition 6.1:} Given $A,B\in M_N$, there exists unitary $U$ such that $B=U^*AU$ iff
\beq{E7.4}{
\trace\bigl(w(A,A^*)\bigr)=\trace\bigl(w(B,B^*)\bigr)}
\end{equation}
for all words $w(\cdot,\cdot)$ of length $\leq\,N^2$.
{\\[.5cm]}
The disparity between (\ref{E7.4}) and (\ref{E7.2}) certainly suggests that $A$ and $B$ might have the same numerical
shadow without being simply related by unitaries. Let us see how this does occur when $N=3$. In [DJ2007], Djokovi\'c
and Johnson refer to a result of Sibirski\v i: the unitary equivalence class of $A\in M_3$ is determined by
\[
\trace(A),\trace(A^2),\trace(AA^*),\trace(A^3),
\]\[
\trace(A^2A^*), \trace\bigl(A^2(A^*)^2\bigr),\mbox{ and }\trace\bigl(A^2(A^*)^2AA^*\bigr),
\]
and this set is minimal. In contrast, (\ref{E7.2}) tells us that the numerical shadow $P_A$ is determined (when
$A\in M_3$) by
\[
\trace(A),\trace(A^2),\trace(AA^*),\trace(A^3),\mbox{ and }\trace(A^2A^*).
\]
Thus we expect to find $A,B\in M_3$ such that $P_A=P_B$ but $A$ and $B$ are not unitarily related.
{\\[.5cm]}
A class of specific examples is provided by
\[
A=\begin{bmatrix}0&x&0\\y&0&0\\z&0&0\end{bmatrix}, B=\begin{bmatrix}0&y&0\\x&0&0\\z&0&0\end{bmatrix}.
\]
Note that $\det(I-(sA+tA^*))=1-st(|x|^2+|y|^2+|z|^2)-s^2xy-t^2\overline{xy}$. Since this expression is symmetric in $x,y$,
(\ref{E7.3}) tells us that $P_A=P_B$. Consider the choice $x=0,y=z=1$: then $A$ has rank 1 while $B$ has rank 2. Clearly
$B$ is not unitarily similar to $A$ or to $A^t$.
{\\[.5cm]}
\textbf{Remark:} The common numerical shadow of these $A,B$ is identified explicitly in section 7, because $B$ is unitarily
equivalent to the Jordan nilpotent $J_3$.

\section{Numerical shadows of Jordan nilpotents $J_N$}

Here we compute explicit shadow densities for certain special
matrices, focusing on the Jordan nilpotent $J_N$, ie $J_N\in M_N(\dC)$ with 1's on the superdiagonal and 0's elsewhere. Of course,
the discussion of the $2\times2$ case in section 2 applies to $J_2$ and shows that the planar density of the shadow $P_{J_2}$
at $z\in\dC$ is $f_2(|z|)$ where
\[
f_2(r)=\frac1{2\pi (1/2)^2 \sqrt{1-4r^2}}=\frac2{\pi\sqrt{1-4r^2}},
\]
since $W(J_2)$ is a disc of radius $1/2$.
The shadow density for $J_3$ can be computed by several methods but here we'll do it as the simplest case of a general
method that exploits the moment techniques from section 5. We shall see that the shadow density for $J_N$ is an alternating
sum of densities supported on discs with centre at 0 and with various radii, the largest being $\cos(\pi/(N+1))$. This is a striking
development beyond the well--known numerical radius formula: $w(J_N)=\cos(\pi/(N+1))$ (see [DH1988] for information about the
numerical radii of certain matrices with simple structure; for more, see [HS2010]).
{\\[.5cm]}
Observe first that the shadow measure $P_{J_N}$ is certainly circularly symmetric about 0; in fact $J_N$ and $e^{i\theta}J_N$
are unitarily similar (use $U=\diag(1,e^{i\theta},e^{i2\theta},\dots)$). Thus, from Proposition 5.3, we have
\beq{8--1}{
\sum_{m=0}^\infty \frac{(N)_m}{m!m!}s^mt^m\nu_{mm}(J_N)=\mbox{det}^{-1}(I_N-(sJ_n+tJ_N^*)).
}
\end{equation}
We may take $s=t$ and identify $\nu_{mm}(J_N)$ via the coefficient of $t^{2m}$ in $\mbox{det}^{-1}\bigl(I_N-t(J_N+J_N^*)\bigr)$.
Now the eigenvalues of $J_N+J_N^*$ are well--known:
\[
2\cos\Bigl(\frac{k\pi}{N+1}\Bigr)\quad(k=1,2,\dots,N).
\]
[Some say that this was the first nontrivial eigenvalue problem ever solved, and that it goes all the way back to Cauchy.]
{\\[.5cm]}
Thus the RHS of (\ref{8--1}) (for $s=t$) can be calculated explicitly. The details appear in the proof of the
following proposition.
{\\[.5cm]}
\textbf{Proposition 7.1:} For each $N\geq2$ and $m=0,1,\dots$
\[
\nu_{mm}(J_N)=\sum_{k=1}^{\llcorner N/2\lrcorner}c_k\bigl(\cos^2(\frac{k\pi}{N+1})\bigr)^m\frac{m!}{(\frac N2)_m}\frac{m!}{(\frac{N+1}2)_m},
\]
where
\[
c_k=(-1)^{k-1}\frac{2^{N+1}}{N+1}\sin^2\bigl(\frac{k\pi}{N+1}\bigr)
\Bigl(\cos\bigl(\frac{k\pi}{N+1}\bigr)\Bigr)^{N-1}.
\]
\textbf{Proof:} By (\ref{8--1}), with (small) $s=t$,
\[
\sum_{m=0}^\infty \frac{(N)_m}{m!m!}t^{2m}\nu_{mm}(J_N)=\frac1{\prod_{k=1}^N
\bigl(1-2t\cos(\frac{k\pi}{N+1})\bigr)}.
\]
Now the $\cos(\frac{k\pi}{N+1})$ are the roots of the monic polynomial $C_N(x)$, where
\[
C_N(\cos\theta)=\frac1{2^N}\frac{\sin(N+1)\theta}{\sin\theta}
\]
(a version of the Chebyshev polynomials of the second kind). Thus
\[
\frac1{C_N(x)}=\frac1{\prod_{k=1}^N\bigl(x-\cos(\frac{k\pi}{N+1})\bigr)}=
\sum_{k=1}^N\frac{a_k}{\bigl(x-\cos(\frac{k\pi}{N+1})\bigr)},
\]
where the coefficients $a_k$ in the partial fraction decomposition are given by $a_k=1/C_N'(\cos(\frac{k\pi}{N+1}))$.
Using the formula for $C_N(\cos\theta)$ we find that
\[
a_k=(-1)^{k-1}\frac{2^N\sin^2(\frac{k\pi}{N+1})}{N+1}.
\]
We now have
\begin{eqnarray*}
\sum_{m=0}^\infty \frac{(N)_m}{m!m!}t^{2m}\nu_{mm}(J_N) & = & \frac1{(2t)^N\prod_{k=1}^N((1/2t)-\cos(\frac{k\pi}{N+1}))} = \frac1{(2t)^N}\frac1{C_N(1/2t)}\\ 
&=&\frac1{t^N}\sum_{k=1}^N\frac{(-1)^{k-1}}{N+1}
\sin^2(\frac{k\pi}{N+1})\frac1{\bigl((1/2t)-\cos(\frac{k\pi}{N+1})\bigr)}\\
&=&\frac1{t^{N-1}}\frac2{N+1}\sum_{k=1}^N(-1)^{k-1}\sin^2(\frac{k\pi}{N+1})\sum_{j=0}^\infty\bigl(2t\cos(\frac{k\pi}{N+1})\bigr)^j.
\end{eqnarray*}
Evidently the summed coefficients for $j$ odd and for $j<N-1$ are 0 [we need not worry about how this happens!] so
that the term in $t^{2m}$ for the final expression corresponds to $j=2m+N-1$. Thus
\[
\nu_{mm}(J_N)=\sum_{k=1}^N(-1)^{k-1}\frac{2^N}{N+1}
\sin^2\bigl(\frac{k\pi}{N+1}\bigr)
\Bigl(\cos\bigl(\frac{k\pi}{N+1}\bigr)\Bigr)^{N-1}(4\cos^2(\frac{k\pi}{N+1}))^m
\frac{m!m!}{(N)_{2m}}.
\]
Note that the terms for $k$ and $N-k+1$ are the same and that $(N)_m$ may be replaced by $(\frac{N}2)_m(\frac{N+1}2)_m2^{2m}$.
Then $\nu_{mm}(J_N)$ is given by
\[
\sum_{k=1}^{\llcorner N/2\lrcorner}(-1)^{k-1}\frac{2^{N+1}}{N+1}
\sin^2\bigl(\frac{k\pi}{N+1}\bigr)\bigl(\cos(\frac{k\pi}{N+1})\bigr)^{N-1}
\bigl(\cos^2(\frac{k\pi}{N+1})\bigr)^m
\frac{m!m!}{(\frac{N}2)_m(\frac{N+1}2)_m},
\]
(the additional factor of 2 is correct even if $N$ is odd because then $\cos(\frac{k\pi}{N+1})=0$ for $k=(N+1)/2$).
QED
{\\[.5cm]}
The value of Proposition 7.1 lies in the possibility of identifying explicitly those densities with moments
\[
b^m\frac{m!}{(\frac{N}2)_m}\frac{m!}{(\frac{N+1}2)_m}.
\]
\textbf{Proposition 7.2:} Suppose $f(x)$ and $g(x)$ are probability densities on $[0,1]$ with moments
\[
\int_0^1x^mf(x)\,dx=a_m,\quad \int_0^1x^mg(x)\,dx=b_m.
\]
Then\\
(i) for any $b>0$, $(1/b)f(x/b)$ is a probability density on $[0,b]$ with moments $b^ma_m$,\\
and\\
(ii) a probability density on $[0,1]$ with moments $a_mb_m$ is given by
\[
h(x)=\int_x^1f(s)g(\frac{x}{s})\,\frac{ds}s.
\]
\textbf{Proof:} (i) With the substitution $y=x/b$,
\[
\int_0^bx^m\frac1b f(\frac{x}b)\,dx=\int_0^1b^my^mf(y)\,dy=b^ma_m.
\]
(ii) Consider independent random variables $X,Y$ with $f,g$ as probability densities. Then $XY$ has moments
\[
E((XY)^m)=E(X^m)E(Y^m)=a_mb_m.
\]
Since $X,Y$ have joint density $f(x)g(y)$, Prob$\{XY\leq t\}$ is
\[
\int_0^1f(x)\Bigl(\int_0^{t/x}g(y)\,dy\Bigr)\,dx=\int_0^tf(x)
\Bigl(\int_0^1g(y)\,dy\Bigr)\,dx+\int_t^1f(x)\Bigl(\int_0^{t/x}g(y)\,dy\Bigr)\,dx,
\]
for $t\in[0,1]$. Differentiate with respect to $t$ to obtain the density $h(t)$ for $XY$:
\[
f(t)-f(t)+\int_t^1f(x)g(\frac tx)\,\frac{dx}x.
\]
QED
{\\[.5cm]}
We can now compute the density $F_N(x)$ on $[0,1]$ having moments
\[
\int_0^1x^mF_N(x)\,dx=\frac{m!}{(\frac{N+1}2)_m}\frac{m!}{(\frac{N}2)_m},
\]
for $m=0,1,2,\dots$ and $N\geq2$. For any $\beta>0$ we have the beta--integrals
\[
\int_0^1\beta x^m(1-x)^{\beta-1}\,dx=\frac{m!}{(\beta+1)_m},
\]
for $m=0,1,2,\dots$ (use induction on $m$ via integration by parts).
For $N=2$ take $\beta=1/2$ to see that
\[
F_2(x)=\frac1{2\sqrt{1-x}}.
\]
For $N\geq3$ we apply Proposition 7.2(ii): $F_N(x)=h(x)$ computed with
\[
f(x)=\frac{N-1}2 (1-x)^{\frac{N-3}2},\quad g(x)=\frac{N-2}2 (1-x)^{\frac{N-4}2}.
\]
Consider even $N=2\ell$:
\begin{eqnarray*}
F_{2\ell}(x)&=&\int_x^1\frac{2\ell-1}2(1-s)^{\frac{2\ell-3}2}\frac{2\ell-2}2(1-\frac xs)^\frac{2\ell-4}2 \,\frac{ds}s \\
&=&\frac{(\ell-1)(2\ell-1)}2\int_x^1(1-s)^{\ell-\frac32}(s-x)^{\ell-2}s^{1-\ell}\,ds.
\end{eqnarray*}

With the substitution $s=1-u^2$ we obtain
\[
F_{2\ell}(x)=(\ell-1)(2\ell-1)\int_0^{\sqrt{1-x}} u^{2\ell-2}(1-x-u^2)^{\ell-2}\,\frac{du}{(1-u^2)^{\ell-1}}.
\]
For odd $N=2\ell+1$ we reverse the roles of $f$ and $g$ to obtain
\[
F_{2\ell+1}(x)=\int_x^1\frac{2\ell-1}2(1-s)^{\ell-\frac32}\,\ell 
\bigl(1-\frac xs\bigr)^{\ell-1}\,\frac{ds}s
\]
\[
=\ell(2\ell-1)\int_0^{\sqrt{1-x}} u^{2\ell-2}(1-x-u^2)^{\ell-1}\,\frac{du}{(1-u^2)^\ell}.
\]
The integrals representing $F_N(x)$ are elementary in the sense that they may in principle be computed explicitly (using
partial fractions, for example). In particular,
\[
F_3(x)=\int_0^{\sqrt{1-x}}\, \frac{du}{1-u^2}=\log\frac{1+\sqrt{1-x}}{\sqrt x},
\]
and
\[
F_4(x)=3\int_0^{\sqrt{1-x}} u^2\,\frac{du}{1-u^2}=3\log\frac{1+\sqrt{1-x}}{\sqrt x} -3\sqrt{1-x}.
\]
In fact, there is a recurrence relation for the $F_N(x)$ that makes the calculation of $F_N$ for $N>4$ a simple task;
such matters are discussed at the end of this section.
{\\[.5cm]}
Returning to shadow densities, let the planar density of $P_{J_N}$ at $z\in\dC$ be denoted by $f_N(|z|)$ so that
\[
\nu_{mm}(J_N)=\int_0^{2\pi}\int_o^{w(J_N)} r^{2m}f_N(r)r\,dr\,d\theta=2\pi\int_0^{w(J_N)}r^{2m+1}f_N(r)\,dr.
\]
With the substitution $x=r^2$ we have
\[
\nu_{mm}(J_N)=\pi\int_0^{w^2(J_N)} x^mf_N(\sqrt x)\,dx.
\]
In view of Proposition 7.1 and Proposition 7.2(i),
\[
\pi\int_0^{w^2(J_N)} x^mf_N(\sqrt x)\,dx=
\]\[
\sum_{k=1}^{\llcorner N/2\lrcorner} c_k\int_0^{\cos^2\frac{k\pi}{N+1}} \frac{x^m}{\cos^2\frac{k\pi}{N+1}}
 F_N\bigl(\frac{x}{\cos^2\frac{k\pi}{N+1}}\bigr)\,dx,
\]
where $c_k$ are as in Proposition 7.1. Thus $\pi f_N(\sqrt x)$ and
\[
\sum_{k=1}^{\llcorner N/2\lrcorner}  \frac{c_k}{\cos^2\frac{k\pi}{N+1}} F_N\bigl(\frac{x}{\cos^2\frac{k\pi}{N+1}}\bigr)
\]
coincide (since they have the same moments). We obtain
\[
f_N(x)=\frac1\pi \sum_{k=1}^{\llcorner N/2\lrcorner}  \frac{c_k}{\cos^2\frac{k\pi}{N+1}} F_N\bigl(\frac{x^2}{\cos^2\frac{k\pi}{N+1}}\bigr),
\]
an (alternating) sum of densities supported on $[0,\cos\frac{k\pi}{N+1}]$; in particular we have a greatly refined version of
the result $w(J_N)=\cos\frac{\pi}{N+1}$ ($=\max_k\cos\frac{k\pi}{N+1}$).
{\\[.5cm]}
In summary, we have proved
{\\[.5cm]}
\textbf{Proposition 7.3:} The radial density $f_N(r)$ of $P_{J_N}$ is given by
\[
f_N(r)=\frac1\pi\sum_{k=1}^{\llcorner N/2\lrcorner} (-1)^{k-1}\frac{2^{N+1}}{N+1}\sin^2\frac{k\pi}{N+1}\bigl(\cos\frac{k\pi}{N+1}\bigr)^{N-3}
F_N\bigl(\frac{r^2}{\cos^2\frac{k\pi}{N+1}}\bigr)
\]
for any $N\geq2$.
{\\[.5cm]}
For $N=2$ we see again that
\[
f_2(r)\quad=\Bigl(\frac8{3\pi}\sin^2\frac\pi3\cos^{-1}\frac\pi3\frac1{2\sqrt{1-r^2/\cos^2\frac\pi3}}\Bigr)=\frac2\pi\frac1{\sqrt{1-4r^2}}.
\]
Likewise
\[
f_3(r)=\frac4\pi\sin^2\frac\pi4\,F_3 \bigl(\frac{r^2}{\cos^2\frac\pi4}\bigr)=\frac2\pi\log\frac{1+\sqrt{1-2r^2}}{\sqrt2r}.
\]
For $N>3$ the radial density combines densities on discs of several different radii. For example,
\[
f_4(r)=\frac{32}{5\pi}\Bigl(\sin^2\frac\pi5\,\cos^2\frac\pi5\,
  F_4\bigl(\frac{r^2}{\cos^2\frac\pi5}\bigr)
-\sin^2\frac{2\pi}5\,\cos^2\frac{2\pi}5\,
F_4\bigl(\frac{r^2}{\cos^2\frac{2\pi}5}\bigr)\Bigr)
\]\[
=\frac1{5\pi}\Bigl((5+\sqrt5)F_4\bigl(\frac{8r^2}{3+\sqrt5}\bigr)-
(5-\sqrt5)F_4\bigl(\frac{8r^2}{3-\sqrt5}\bigr)\Bigr).
\]
We shall see that the functions $F_N$ are the basic building blocks for many circularly symmetric numerical shadows.
Hence it will be worthwhile to explore their properties more thoroughly. To this end, we introduce the following
hypergeometric series:
\[
H(a,b;c;t)=\sum_{j=0}^\infty \frac{(a)_j(b)_j}{(c)_j j!}t^j;
\]
in this context $H$ is often denoted by $_2F_1$.
Here we may assume that the parameters $a,b,c$ are real and that $c\neq 0,-1,-2,\dots$. Note that the series converges
absolutely for $|t|<1$ since it has the form $\sum d_jt^j$ where
\[
|d_{j+1}/d_j|=|(a+j)(b+j)/((c+j)(j+1))|\to_j1.
\]
Recall the Gauss summation formula, which tells us that the series converges also for $t=1$ whenever $a,b\geq0$ and
$c-a-b>0$ and that, in such a case,
\[
H(a,b;c;1)=\frac{\Gamma(c-a-b)\Gamma(c)}{\Gamma(c-a)\Gamma(c-b)}.
\]
Given $\beta\geq0,\delta>0$ and $k=1,2,\dots$, let
\beq{E8.2}{
G(x)=\frac{\Gamma(\beta+1)\Gamma(\delta+k)}{\Gamma(\beta+\delta)(k-1)!} (1-x)^{\beta+\delta-1}H(\delta,\beta+1-k;\beta+\delta;1-x),}
\end{equation}
for $0<x<1$.
{\\[.5cm]}
\textbf{Proposition 7.4:} The function defined by (\ref{E8.2}) is a probability density on $(0,1)$ with moments
\[
\int_0^1 x^m G(x)\,dx=\frac{m!\;(k)_m}{(\beta+1)_m(\delta+k)_m}\quad(m=0,1,2,\dots).
\]
\textbf{Proof:} Evaluating the beta--functions $\int_0^1x^m(-x)^{\beta+\delta-1+j}\,dx$ as
\[
\frac{m!\;\Gamma(\beta+\delta+j)}{\Gamma(\beta+\delta+j+m+1)}=\frac{m!\;\Gamma(\beta+\delta)(\beta+\delta)_j}
{\Gamma(\beta+\delta)(\beta+\delta)_{m+1+j}}
\]\[
=\frac{m!\;(\beta+\delta)_j}{(\beta+\delta)_{m+1}(\beta+\delta+m+1)_j},
\]
we see that
\[
\int_0^1 x^m G(x)\,dx=\frac{m!\;\Gamma(\beta+1)\Gamma(\delta+k)}
{\Gamma(\beta+\delta)(k-1)!\;(\beta+\delta)_{m+1}}
\sum_{j=0}^\infty\frac{(\delta)_j(\beta+1-k)_j}{(\beta+\delta+m+1)_j\;j!}.
\]
Using the Gauss summation formula, we obtain
\begin{eqnarray*}
&& \frac{m!\;\Gamma(\beta+1)\Gamma(\delta+k)}{\Gamma(\beta+\delta)(k-1)!\; (\beta+\delta)_{m+1}}\cdot
\frac{\Gamma(m+k)\Gamma(\beta+\delta+m+1)}{\Gamma(\beta+1+m)\Gamma(\delta+k+m)} \\
&& =\frac{m!\;\Gamma(\beta+1)\Gamma(\delta+k)}{\Gamma(\beta+\delta+m+1)(k-1)!}
\cdot \frac{(k-1+m)!\;\Gamma(\beta+\delta+m+1)}
{\Gamma(\beta+1)(\beta+1)_m\Gamma(\delta+k)(\delta+k)_m}\\
&& =\frac{m!\;(k)_m}{(\beta+1)_m(\delta+k)_m}.\quad\mbox{QED}
\end{eqnarray*}
Taking $k=1, \beta=(N-2)/2,\delta=(N-1)/2$ we see in particular that
\beq{E8.3}{
F_N(x)=\frac{\Gamma(\frac{N}2)\Gamma(\frac{N+1}2)}
{\Gamma(N-\frac32)}(1-x)^{N-\frac52}H(\frac{N-1}2,\frac{N-2}2;N-\frac32;1-x).}
\end{equation}
{\\[.5cm]}
Several useful recurrence relations will follow from the following general recurrence for $H$.
{\\[.5cm]}
\textbf{Lemma 7.5:} If $a,b$ are real and $b>1$ then
\[
H(a-\frac12,a-1;b;t)-H(a-1,a-\frac32;b-1;t)
\]\[
=\frac{(2b+1-2a)(a-1)}{2b(b-1)}\cdot tH(a,a-\frac12;b+1;t).
\]
This lemma may be verified by a careful comparison of the terms involving $t^{j+1}$ ($j=0,1,\dots$).
{\\[.5cm]}
Using (\ref{E8.3}) and invoking the lemma with $t=1-x,a=(N+1)/2,b=N-\frac12$, we obtain the recurrence relation for
$F_{N+2}(x)$ ($N\geq2$):
\beq{E8.4}{
F_{N+2}(x)=\frac{N+1}{(N-1)^2}\big((2N-3)F_{N+1}(x)-(1-x)N\,F_N(x)\big).}
\end{equation}
In fact, then, each $F_N(X)$ has the form $a_N(x)F_2(x)+b_N(x)F_3(x)$ for certain polynomials $a_N,b_N$.
{\\[.5cm]}
We extend the definition of $F_N(x)$ by setting it equal to 0 for $x\geq1$; this is the natural continuous extension
except that
$F_2(x)=(2\sqrt{1-x})^{-1}\uparrow\infty$ as $x\uparrow1$. From (\ref{E8.3}) it is clear that for $N\geq3$
we  have $F_N(x)\to 0$ as $x\uparrow1$; hence $a_N(1)=0$ for $N\geq3$. We may also examine the behavior of
the functions $F_N$ at 0: $F_3(x)=\log(1+\sqrt{1-x})-\frac12\log x\uparrow\infty$ as $x\downarrow0$, whereas
$F_2(x)\to\frac12$. On the other hand (\ref{E8.3}) shows that, for $N\geq3$, $F_N(x)$ tends to a constant
times $H\bigl((N-1)/2,(N-2)/2;N-3/2;1\bigr)$ as $x\downarrow0$; the Gauss summation formula tells us that this limit
is $+\infty$ (since $c-a-b=0$). Thus $b_N(0)>0$ for $N\geq3$, and $F_N(x)$ grows like $-\log x$ as $x\downarrow0$

\section{Rotation--invariant shadows}

Here we shall see that the methods of section 7 extend to determine
explicit densities for all rotation--invariant numerical shadows. These are shadows of $A\in M_N$ such that
$A$ and $e^{i\theta}A$ have the same shadow for all real $\theta$. Characterizing such $A$ in terms of moments
is easy: $\nu_{jk}(A)=0$ whenever $j\neq k$. More elusive are characterizations directly in terms of $A$.
{\\[.5cm]}
Simple examples are provided by the ``superdiagonal'' matrices: ie $A$ such that $a_{ij}=0$ unless $j=i+1$. For
such $A$ we actually have $e^{i\theta}A$ unitarily similar to $A$: let $U=\diag(1,e^{i\theta},e^{i2\theta},\dots)$;
then $U^*AU=e^{i\theta}A$. The Jordan nilpotents $J_N$ are special cases of these superdiagonal matrices.
{\\[.5cm]}
More generally, consider the incidence graph $G(A)$ of $A\in M_N$: vertices are $\{1,2,\dots,N\}$ and $i,j$ are joined
by an edge iff $a_{ij}\neq 0$. The interesting case in this context is when $G(A)$ consists of disjoint chains (no cycles
are allowed; in particular, $A$ has zero diagonal). One can see that this condition is equivalent to requiring that $A$
have zero diagonal, have no more than two nonzero entries in each cross--shaped region formed by the $k$--th row and
the $k$--th column, and that $G(A)$ have no cycles.
{\\[.5cm]}
\textbf{Proposition 8.1:} If $G(A)$ consists of disjoint chains, then $A$ and $e^{i\theta}A$ are unitarily similar (so that
$A$ has rotation--invariant shadow).\\
\textbf{Proof:} Consider the unitary $U=\diag(u)$ where $u_{j_k}=e^{ik\theta}$ for each chain
\[
j_1\to j_2\to\dots \to j_k\to\dots\to j_K
\]
of $G(A)$ (it does not matter which orientation of the chain is chosen). Set $u_j=1$ for any $j$ that does not occur in any
of the chains that make up $G(A)$. Note that
\[
(U^*AU)_{j_k,j_{k+1}} = e^{ik\theta}a_{j_k,j_{k+1}}e^{i(k+1)\theta}=e^{i\theta}a_{j_k,j_{k+1}}.
\]
Since other entries of $A$ are 0, we do have $U^*AU=e^{i\theta}A$. QED
{\\[.5cm]}
This proposition applies, for example, to superdiagonal $A$ as well as to strictly upper--triangular $A$ that
are ``subpermutation'' matrices, ie have at most one nonzero entry in each row and in each column.
{\\[.5cm]}
The next proposition notes that $A$ with rotation--invariant shadow must be nilpotent, so that it is unitarily
similar to a strictly upper--triangular matrix (Schur form).
{\\[.5cm]}
\textbf{Proposition 8.2:} If $A\in M_N$ has rotation--invariant numerical shadow, then all eigenvalues are 0.\\
\textbf{Proof:} Putting $t=0$ in (\ref{E6.5}), we see that $h_n(\ld(A))=0$ ($n\geq1$); indeed, this is  the case whenever
$\nu_{n,0}(A)=0$. From (\ref{E6.6}) and (\ref{E6.7}) we conclude that $p_n(\ld(A))=0$ for $n\geq1$.
Thus $\sum_1^N p(\ld_k)=0$ for any polynomial $p(x)$ with $p(0)=0$. Suppose $\ld_i$ occurs with multiplicity $m$.
Let
\[
p(x)=x\prod_{\ld_j\neq\ld_i}(x-\ld_j);
\]
then
\[
0=\sum_1^N p(\ld_k)=m\ld_i\prod_{\ld_j\neq\ld_i}(\ld_i-\ld_j),
\]
so that $\ld_i=0$. QED
{\\[.5cm]}
\textbf{Proposition 8.3:} The matrix $A\in M_N$ has rotation--invariant numerical shadow\\
\textbf{iff}
\[
\mbox{(i)}\quad\det\bigl(I-(sA+tA^*)\bigr)\mbox{ is a function of }st;
\]
\textbf{iff}
\[
\mbox{(ii) $A$ is nilpotent and}\quad\trace\bigl(P_{k,j}(A,A^*)\bigr)=0\quad(\frac{k}2<j<k\leq N).
\]
\textbf{Proof:} In view of (\ref{E6.8}), (i) is equivalent to $\nu_{k,j}(A)=0$ for $k\neq j$. To see that
(ii) follows from rotation--invariance, invoke Proposition 8.2 and apply (\ref{E7.2}) with $B=e^{i\theta}A$
to obtain
\beq{E9.1}{
\trace\bigl (P_{k,j}(A,A^*)\bigr)=e^{i[j-(k-j)]\theta}\trace\bigl(P_{k,j}(A,A^*)\bigr).}
\end{equation}
When $j\neq k/2$, this cannot hold (for all $\theta$) unless
 $\trace\bigl(P_{k,j}(A,A^*)\bigr)=0$. For the converse,
note that (\ref{E9.1}) is automatic when $j=k/2$ and that nilpotence ensures that both sides of (\ref{E7.2})
are zero also when $j=k$. For $j<k/2$, note that
\[
\trace(P_{k,j}(A,A^*))=\trace\bigl((P_{k,k-j}(A,A^*)\bigr)^*)=\overline{\trace(P_{k,k-j}(A,A^*))}.
\]
QED
{\\[.5cm]}
When $N=3$, either (i) or (ii) easily implies that the upper--triangular form of $A\in M_3$ with rotation--invariant
shadow is
\[
\begin{bmatrix} 0&x&y\\0&0&z\\0&0&0\end{bmatrix},
\]
where at least one of $x,y,z$ is zero. For example, the only condition in (ii) is that $\trace(P_{3,2}(A,A))=0$, ie
that $\trace(A^2A^*)=0$, and one easily computes $\trace(A^2A^*)=x\overline{y}z$. Note that $A$ and $B=e^{i\theta}A$ are
unitarily similar. One can appeal to Proposition 8.1 to see this or use Sibirski\v i's list of words (mentioned in section 6):
all the traces are automatically the same for nilpotent $A$ and $B$ except that $\trace(A^2A^*)=\trace(B^2B^*)$ requires
$\trace(A^2A^*)=0$.
{\\[.5cm]}
For $N=4$ it is perhaps more convenient to use (ii) to identify those $A$ having rotation--invariant shadow. Let
the upper--triangular form of $A$ be
\[
\begin{bmatrix} 0&a&b&c\\0&0&x&y\\0&0&0&z\\0&0&0&0\end{bmatrix}.
\]
The only conditions in (ii) when $N=4$ are $\trace(P_{3,2}(A,A^*))=0$ and $\trace(P_{4,3}(A,A^*))=0$,
ie $\trace(A^2A^*)=0$ and $\trace(A^3A^*)=0$. Computing these traces we find that $A$ has rotation--invariant
shadow iff
\beq{E9.2}{
ax\overline{b}+(ay+bz)\overline{c}+x\overline{y}z=0 \mbox{ and } ax^2z\overline{c}=0.}
\end{equation}
{\\[.5cm]}
\textbf{Remark 8.4:} Although the earlier examples of $A$ with rotation--invariant shadow were also unitarily
similar to $e^{i\theta}A$, the analysis (above) of the $4\times4$ case shows that this is not necessary. If
$A$ and $e^{i\theta}A$ are unitarily similar we must have $\trace(A^3(A^*)^2)=0$, ie $ax^2z\overline{(ay+bz)}=0$,
and this does not follow from (\ref{E9.2}) (eg take $c=0$, $a=b=x=1$, $y=-1/2$, and $z=2$).
{\\[.5cm]}
If $A\in M_N$ has rotation--invariant shadow, the relation (\ref{E6.8}) simplifies:
\beq{E9.3}{
\sum_{m=0}^\infty t^{2m}\frac{(N)_{2m}}{m!\; m!}\nu_{mm}(A)=\mbox{det}^{-1}(I-2t\mbox{Re}A)}
\end{equation}
(for all sufficiently small real $t$), where Re$A$ is the Hermitian $(A+A^*)/2$.
If $\ld_1,\dots,\ld_K$ are the nonzero eigenvalues (real) of Re$A$, the RHS of (\ref{E9.3}) is
$\bigl(\prod_{k=1}^K (1-2t\ld_k)\bigr)^{-1}$; since the LHS is a function of $t^2$, these eigenvalues come in $\pm$ pairs.
We may assume that $\ld_1,\dots,\ld_p$ are the positive eigenvalues of Re$A$ so that the spectrum of Re$A$ is
\[
(\ld_1,\ld_2,\dots,\ld_p,0,-\ld_1,-\ld_2,\dots,-\ld_p),
\]
where 0 has multiplicity $N-2p$. Note that $p\geq1$ unless $A=0_N$, since Re$A=0$ implies that $A$ is skew--Hermitian
and Proposition 8.2 then implies that $A=0$. We may therefore write (\ref{E9.3}) in the following form:
\beq{E9.4}{
\sum_{m=0}^\infty t^{2m}\frac{(N)_{2m}}{m!\; m!}\nu_{mm}(A)=\frac1{\prod_{j=1}^p(1-4t^2\ld_j^2)}.}
\end{equation}
The methods of section 7 extend most readily to the case where $\ld_1,\dots,\ld_p$ are distinct (as they are
for $A=J_N$, where $p=\llcorner\frac{N}2\lrcorner$ and $\ld_j=\cos(\frac{j\pi}{N+1})$). The following more
general proposition replaces Proposition 7.3.
{\\[.5cm]}
\textbf{Proposition 8.5:} If $0\neq A\in M_N$ has rotation--invariant shadow and the positive eigenvalues of
Re$A$ are the distinct $\ld_1,\dots,\ld_p$ then the planar shadow density at each $z$ with $|z|=r$ is given by
\beq{E9.5}{
f(r)=\frac1\pi \sum_{k=1}^p \frac{\ld_k^{2(p-2)}}{\prod_{1\leq j\leq p,j\neq k}(\ld_k^2-\ld_j^2)} F_N\bigl(\frac{r^2}{\ld_k^2}\bigr),}
\end{equation}
where the function $F_N$ is computable as in section 7.
{\\[.5cm]}
\textbf{Remark 8.6:} To see that Proposition 7.3 is a special case of Proposition 8.5, recall from the proof of
Proposition 7.1 that
\[
(-1)^{k-1}\frac{2^N}{N+1} \sin^2(\frac{k\pi}{N+1})=\frac1{C_N'\bigl(\cos(\frac{k\pi}{N+1})\bigr)}
\]
where $C_N(x)=\prod_{j=1}^N\bigl(x-\cos(\frac{j\pi}{N+1})\bigr)$. Thus
\beq{E9.6}{
(-1)^{k-1}\frac{2^N}{N+1} \sin^2\bigl(\frac{k\pi}{N+1}\bigr)=\frac1{\prod_{1\leq j\leq N,j\neq k}\bigl(\cos(\frac{k\pi}{N+1})-\cos(\frac{j\pi}{N+1})\bigr)}\,\,.}
\end{equation}
When $\ld_j=\cos(\frac{j\pi}{N+1})$, $j=1,2,\dots,N$, the positive values are $\ld_1,\dots,\ld_p$ with $p=\llcorner\frac{N}2\lrcorner$.
{\\[.5cm]}
Suppose first that $N$ is odd; then $p=(N-1)/2$ and for $k\leq p$ we have
\[
\frac{\ld_k^{2(p-2)}}{\prod_{1\leq j\leq p,j\neq k}(\ld_k^2-\ld_j^2)}=
\frac{\ld_k^{N-5}}{\prod_{1\leq j\leq p,j\neq k}(\ld_k-\ld_j)(\ld_k+\ld_j)}
\]\[
=\frac{\ld_k^{N-5}}{\prod_{1\leq j\leq p,j\neq k}(\ld_k-\ld_j)(\ld_k-\ld_{N+1-j})}
=\frac{\ld_k^{N-5}\cdot(\ld_k-0)\cdot(\ld_k-\ld_{N+1-k})}{\prod_{1\leq j\leq N,j\neq k}(\ld_k-\ld_j)}\,\,,
\]
since $\ld_{\frac{N+1}2}=0$. Thus, in view of (\ref{E9.6}),
\[
\frac{\ld_k^{2(p-2)}}{\prod_{1\leq j\leq p,j\neq k}(\ld_k^2-\ld_j^2)}
=\ld_k^{N-5}\cdot\ld_k\cdot2\ld_k(-1)^{k-1}\frac{2^N}{N+1} \sin^2(\frac{k\pi}{N+1})
\]\[
=(-1)^{k-1}\frac{2^{N+1}}{N+1} \sin^2\bigl(\frac{k\pi}{N+1}\bigr)\ld_k^{N-3}\,\,,
\]
and Proposition 7.3 follows from Proposition 8.5. The argument for even $N$ is similar.
{\\[.5cm]}
\textbf{Proof of Proposition 8.5:} Since $\ld_1^2,\dots,\ld_p^2$ are distinct,
\[
\frac1{\prod_{j=1}^p(x-\ld_j^2)}=\sum_{k=1}^p\frac{b_k}{(x-\ld_k^2)}\,\,,
\]
where
\[
b_k=\frac1{\prod_{j\neq k}(\ld_k^2-\ld_j^2)}\,\,.
\]
With $x=1/4t^2$, the RHS of (\ref{E9.4}) becomes
\[
\frac1{(4t^2)^p}\sum_{k=1}^p\frac{b_k}{(x-\ld_k^2)}\,\,,
\]
which we may write as
\[
\frac1{(4t^2)^{p-1}}\sum_{k=1}^p b_k(1+4t^2\ld_k^2+(4t^2)^2\ld_k^4+\dots)\,\,.\]
Comparing terms involving $t^{2m}$ with (\ref{E9.4}) we see that
\[
\nu_{mm}(A)=\frac{m!m!}{(N)_{2m}}\sum_{k=1}^p b_k 4^m\ld_k^{2(m+p-1)} .
\]
In terms of the radial density $f$ we have
\[
\pi\int_0^\infty x^mf(\sqrt{x})\,dx=\frac{m!m!}{(\frac{N}2)_m(\frac{N+1}2)_m2^{2m}}\sum_{k=1}^p b_k 4^m\ld_k^{2(m+p-1)},
\]
so that (in view of the moments that $F_N$ was designed to have)
\begin{eqnarray*}
\pi\int_0^\infty x^mf(\sqrt{x})\,dx &=& (\int_0^1 x^mF_N(x)\,dx)\sum_{k=1}^p b_k \ld_k^{2(m+p-1)} \\
&=& \sum_{k=1}^p b_k\ld_k^{2(p-1)}((\ld_k^2)^m\int_0^1x^mF_N(x)\,dx).
\end{eqnarray*}
Applying Proposition 7.2(i),
\begin{eqnarray*}
\pi\int_0^\infty x^mf(\sqrt{x})\,dx &=& \sum_{k=1}^p b_k\ld_k^{2(p-1)}\int_0^{\ld_k^2}x^m\frac1{\ld_k^2}F_N(x/\ld_k^2)\,dx \\
&=& \int_0^\infty x^m\Bigl(\sum_{k=1}^p b_k\ld_k^{2(p-2)}F_N(x/\ld_k^2)\Bigr)\,dx.
\end{eqnarray*}
Since all moments coincide,
\[
f(\sqrt{x})=\frac1\pi \sum_{k=1}^p b_k\ld_k^{2(p-2)}F_N(x/\ld_k^2).
\]
QED
{\\[.5cm]}
\textbf{Remark 8.7:} Whether or not $\ld_1,\dots,\ld_p$ are distinct, (\ref{E9.4}) shows that the
shadow measure depends only on the $\ld_k$. Thus $B=\oplus_{k=1}^p 2\ld_k J_2\oplus 0_{N-2p}$
has the same shadow as $A$, because the positive eigenvalues of Re$B$ are also $\ld_1,\dots,\ld_p$.
{\\[.5cm]}
All rotation--invariant numerical shadows are obtained as shadows of the simple superdiagonal matrices
\[
B=\oplus_{k=1}^p 2\ld_k J_2\oplus 0_{N-2p},
\]
where $\ld_1,\dots,\ld_p >0$. For example, $A=J_3$ has the same numerical shadow as
\[
B=\begin{bmatrix}0&\sqrt2&0\\0&0&0\\0&0&0\end{bmatrix}.
\]
Here we have another simple example of a pair of matrices with different ranks but the same shadow
(compare the discussion at the end of section 6).
{\\[.5cm]}
One way to deal with the case of repetitions among $\ld_1,\dots,\ld_p$ is to follow the method of
Proposition 8.5 but with the necessarily more complicated partial fraction decomposition. Suppose the distinct values
are $\mu_1,\dots,\mu_n$ and that $\mu_i$ occurs with multiplicity $k_i$; then $p=\sum_1^n k_i$ and det$^{-1}\bigl(I-t(A+A^*)\bigr)$ is
\beq{E9.7}{
\frac1{\prod_{i=1}^n(1-4t^2\mu_i^2)^{k_i}}
=\sum_{i=1}^n\,\,\sum_{j=0}^{k_i-1}\frac{\alpha_{ij}}{(1-4t^2\mu_i^2)^{k_i-j}},}
\end{equation}
for certain constants $\alpha_{ij}$.
{\\[.5cm]}
\textbf{Remark 8.8:} 
The $\alpha_{ij}$ are functions of the eigenvalue data. Computationally effective expressions for these functions are available: see [Hn1974, pp. 553-562].
{\\[.5cm]}
Let $R_{N,k}(y)$ be defined by
\beq{E9.8}{
\frac{\Gamma(\frac{N}2)\Gamma(\frac{N+1}2)}{\Gamma(N-\frac12-k)(k-1)!}(1-y)^{N-\frac32-k}\;
H\Bigl(\frac{N+1}2-k,\frac{N}2-k;N-\frac12-k;1-y\Bigr),}
\end{equation}
with the understanding that $R_{N,k}(y)=0$ for $y\geq1$. In view of Proposition 7.4,
\beq{E9.9}{
\int_0^1 y^MR_{N,k}(y)\,dy = \frac{m!\;(k)_m}{(\frac{n}2)_m(\frac{N+1}2)_m}.}
\end{equation}
\textbf{Proposition 8.9:} If $0\neq A\in M_n$ has rotation--invariant shadow and the positive eigenvalues of Re$A$ are
distinct $\mu_1,\dots,\mu_n$ where $\mu_i$ has multiplicity $k_i$, then the planar shadow density at each $z$ with $|z|=r$
is given by
\beq{E9.10}{
f(r)=\frac1\pi\sum_{i=1}^n\big(\sum_{j=0}^{k_i-1} \alpha_{ij}\frac1{\mu_i^2} R_{N,k_i-j}(\frac{r^2}{\mu_i^2})\big),}
\end{equation}
where $\alpha_{ij}$ are the constants occurring in (\ref{E9.7}).\\
\textbf{Proof:} From (\ref{E9.4}) we obtain
\[
\sum_{m=0}^\infty t^{2m}\frac{(N)_{2m}}{m!\;m!} \nu_{mm}(A)=\sum_{i=1}^n\Big(\sum_{j=0}^{k_i-1}\frac{\alpha_{ij}}{(1-4t^2\mu_i^2)^{k_i-j}}\Big)=
\]\[
\sum_{i=1}^n\Big(\sum_{j=0}^{k_i-1}\alpha_{ij}\,\,\sum_{m=0}^\infty \frac{(k_i-j)_m}{m!} (4t^2\mu_i^2)^m\Big),
\]
where we have used the binomial theorem to express $(1-4t^2\mu_i^2)^{k_i-j}$ (for small t). Comparing coefficients,
\begin{eqnarray*}
\nu_{mm}(A) &=& \frac{m!}{(N)_{2m}/2^{2m}}\sum_{i=1}^n\Big(\sum_{j=0}^{k_i-1} \alpha_{ij}(k_i-j)_m(\mu_i^2)^m\Big) \\ 
&=& \sum_{i=1}^n\Big(\sum_{j=0}^{k_i-1} \alpha_{ij}\,(\mu_i^2)^m\frac{m!\; (k_i-j)_m}{(\frac{N}2)_m(\frac{N+1}2)_m}\Big).
\end{eqnarray*}
In terms of the radial density $f=f_A$, we have
\[
\pi\int_0^\infty x^mf(\sqrt{x})\,dx=\sum_{i=1}^n\Big(\sum_{j=0}^{k_i-1} \alpha_{ij}\,(\mu_i^2)^m\int_0^1y^m R_{N,k_i-j}(y)\,dy \Big)
\]
(recall (\ref{E9.8}) and (\ref{E9.9})). With the substitutions $x=\mu_i^2y$, the RHS becomes
\[
\sum_{i=1}^n\Big(\sum_{j=0}^{k_i-1} \alpha_{ij}\frac1{\mu_i^2}\,\int_0^{\mu_i^2}x^m R_{N,k_i-j}\bigl(\frac{x}{\mu_i^2}\bigr)\,dx \Big).
\]
Because all moments coincide,
\[
f(\sqrt{x})=\frac1\pi\sum_{i=1}^n\Big(\sum_{j=0}^{k_i-1} \alpha_{ij}\frac1{\mu_i^2}\, R_{N,k_i-j}\bigl(\frac{x}{\mu_i^2}\bigr) \Big),
\]
and (\ref{E9.10}) follows. QED
{\\[.5cm]}
\textbf{Remark 8.10:} For example, if all $\ld_k$ have the same value $\mu$, ie $n=1,\mu_1=\mu,k_1=p$, the model matrix is
\[
A=\oplus_1^p \,\,2\mu J_2\oplus\,0_{N-2p}
\]
and the only nonzero $\alpha_{ij}$ in (\ref{E9.7}) is $\alpha_{1,0}=1$. Thus
\[
f(r)=\frac1\pi\frac1{\mu^2} R_{N,p}\bigl(\frac{r^2}{\mu^2}\bigr).
\]
In particular, if $A=\oplus_1^p \,\,2\mu J_2$ (i.e. $N=2p$) we have radial density
\[
f(r)=\frac1\pi\frac1{\mu^2} R_{2p,p}\bigl(\frac{r^2}{\mu^2}\bigr).
\]
Recalling (\ref{E9.8}), we see that
\[
R_{2p,p}(x)=\frac{\Gamma(p)\Gamma(p+\frac12)}{\Gamma(p-\frac12)(p-1)!}(1-x)^{p-\frac32}
\;H\bigl(\frac12,0;p-\frac12;1-x\bigr)
\]\[
=(p-\frac12)(1-x)^{p-\frac32}=(p-\frac12)(\sqrt{1-x})^{2p-3}.
\]
Thus $\oplus_1^p\,\, 2\mu J_2$ has radial density
\[
\frac1\pi\frac1{\mu^2}(p-\frac12)(1-\frac{r^2}{\mu^2})^{p-\frac32}.
\]
\textbf{Remark 8.11:} There is a useful recurrence relation for the functions $R_{N,k}$.
 Apply Lemma 7.5 with
$a=\frac{N}2-k+1$, $b=N-k-\frac12$, $t=1-x$ to see that for $N>2k\geq2$ we have
\[
R_{N+1,k}(x)=\frac{N}{(N-2)(N-2k)}\big((2N-2k-3)R_{N,k}(x)-(N-1)(1-x)R_{N-1,k}(x)\big).
\]
In using this recurrence relation to compute $R_{N,k}$ one would start with $R_{2k,k}$ and $R_{2k+1,k}$. In Remark 8.10
we saw that $R_{2k,k}(x)=(k-\frac12)(\sqrt{1-x})^{2k-3}$; it may also be shown that
\[
R_{2k+1,k}(x)=\frac{k(\frac12)_k}{2(k-1)!}\big((-x)^{k-1}\log(\frac{1+\sqrt{1-x}}{\sqrt{x}})+p_k(x)\sqrt{1-x}\big),
\]
where $p_k(x)$ is an explicitly computable polynomial of degree $k-2$.
{\\[.5cm]}
Further insight into the case of repetitions among the $\ld_1,\dots,\ld_p$ may be gained by considering the limit of
(\ref{E9.5}) from Proposition 8.5 as some of the (initially distinct) $\ld_k$ coalesce. This procedure is legitimate in view
of  the models $A=\oplus_1^p\,\,2\ld_k J_2$, which always have rotation--invariant shadows. In this approach the theory
of divided differences plays an important role. Recall that, given a function $g:[a,b]\to\dR$ and distinct
$y_1,\dots,y_p\in[a,b]$, the divided difference
\[
g[y_1,\dots,y_p]=\sum_{k=1}^p \frac{g(y_k)}{\prod_{j\neq k}(y_k-y_j)}.
\]
We shall appeal to the following facts about such divided differences (compare chapter 4 of [CK1985]):
{\\[.5cm]}
$g[y_1,\dots,y_p]$ is invariant under permutations of the $y_k$;
\beq{E9.11}{
g[y_1,\dots,y_p]=\frac{g[y_1,\dots,y_{p-1}]-g[y_2,\dots,y_p]}{y_1-y_p};}
\end{equation}
if $g$ is $p-1$ times continuously differentiable on $[a,b]$ then
\beq{E9.12}{
\lim \{g[y_1,\dots,y_p]: \mbox{ all } y_k\to y_0\}=\frac{g^{(p-1)}(y_0)}{(p-1)!}.}
\end{equation}
Setting $y_j=1/\ld_j^2$ in Proposition 8.5, we find that
\[
f(r)=\frac1\pi(-1)^{p-1} (\prod_1^p y_j)\,\,g[y_1,\dots,y_p],
\]
where $g(y)=F_N(r^2y)$. In this approach we see that if the distinct positive eigenvalues of Re$A$ are
$\mu_1,\dots,\mu_n$ with multiplicities $k_1,\dots,k_n$, then the radial density $f(r)$ for the numerical shadow
of $A$ may be computed as
\[
f(r)=\frac{(-1)^{p-1}}\pi(\prod_1^n\mu_i^{2k_i})^{-1}\,\,L(k_1,\dots,k_n),
\]
where $p=\sum_1^n k_i$,
\beq{E9.13}{
L(k_1,\dots,k_n)=\lim\{g[y_1,\dots,y_p]:y_j\to 1/\mu_i^2 \mbox{ for }j\in J_i\},}
\end{equation}
and the $J_i$ partition $\{1,2,\dots,p\}$ with $\#(J_i)=k_i$.
{\\[.5cm]}
The relations (\ref{E9.11}) and (\ref{E9.12}) provide us with a sort of $L$--calculus; for example,
\[
L(k,0,\dots,0)=\frac{g^{(k-1)}(1/\mu_1^2)}{(k-1)!}=\frac{r^{2(k-1)}F_N^{(k-1)}(r^2/\mu_1^2)}{(k-1)!},
\]
and
\[
L(k_1,\dots,k_n)=\frac{L(k_1,\dots,k_n-1)-L(k_1-1,\dots,k_n)}{\frac1{\mu_n^2}-\frac1{\mu_n^2}}.
\]
Using those relations repeatedly, we find that
\beq{E9.14}{
L(k_1,\dots,k_n)=\sum_{i=1}^n\Big(\sum_{j=0}^{k_i-1}\frac{\beta_{ij}r^{2j}}{j!} F_N^{(j)}(r^2/\mu_i^2)\Big)}
\end{equation}
for certain constants $\beta_{ij}$. Again (compare Remark 8.8), the $\beta_{ij}$ are functions of the eigenvalue data.
{\\[.5cm]}
Summarizing, we have the following alternate method of computing $f_A(r)$.
{\\[.5cm]}
\textbf{Proposition 8.12:} If $0\neq A\in M_n$ has rotation--invariant shadow and the positive eigenvalues of Re$A$ are
distinct $\mu_1,\dots,\mu_n$ where $\mu_i$ has multiplicity $k_i$, then the planar shadow density at each $z$ with $|z|=r$
is given by
\beq{E9.15}{
f(r)=\frac{(-1)^{p-1}}\pi (\prod_1^n\mu_i^{2k_i})^{-1}\sum_{i=1}^n\big(\sum_{j=0}^{k_i-1}\frac{\beta_{ij}r^{2j}}{j!} F_N^{(j)}(r^2/\mu_i^2)\big),}
\end{equation}
where $p=\sum_1^n k_i$ and $\beta_{ij}$ are the constants found in (\ref{E9.14}).
{\\[.5cm]}
\textbf{Remark 8.13:} A comparison of Propositions 8.9 and 8.12 suggests a relation between $R_{N,k}$ and the derivatives $F_N^{(j)}$.
Indeed, if
\[
A=\oplus_1^p\,\,2J_2\oplus0_{N-2p}
\]
we have $n=1$, $k_1=p$, $\mu_1=1$, $\alpha_{1,0}=1$, and
\[
L(p)=\frac{r^{2(p-1)}}{(p-1)!}\; F_N^{(p-1)}(r^2),
\]
i.e. $\beta_{1,p-1}=1$ and all other $\beta_{ij}=0$. The two forms for the radial density $f_A$ tell us that
\[
R_{N,p}(r^2)=(-1)^{p-1}\frac{r^{2(p-1)}}{(p-1)!} \;F_N^{(p-1)}(r^2),
\]
i.e.
\beq{E9.16}{
R_{N,k}(x)=(-1)^{k-1}\frac{x^{k-1}}{(k-1)!}\; F_N^{(k-1)}(x).}
\end{equation}
This relation between the $R_{N,k}$ and the derivatives of $F_N$ ($=R_{N,1}$) may also be obtained directly by using the identity
\[
H(a,b;c;t)=(1-t)^{c-a-b}H(c-a,c-b;c;t).
\]
\textbf{Remark 8.14:} A study of the behavior of the radial density $f_A(x)$ (when $0\neq A\in M_N$ has rotation--invariant shadow) near
$x=0$ reveals that it has dominant singularity $x^{p-1}\log x$ there (recall that $p$ is the number of positive eigenvalues of Re$A$,
counted with multiplicity) unless $N=2p$, in which case $f_A(x)$ is analytic near $x=0$.

\section{Numerical shadows via the Cartesian decomposition}

In this section we discuss aspects of the numerical shadow related to the
so--called Cartesian decomposition of a matrix $A$ into its Hermitian
components Re$A$ and Im$A={\rm Re}(-iA)$. For example, we investigate the
shadow of a possibly nonnormal matrix by
means of projections onto lines in $\dC$. These projections have
interpretations as shadows of Hermitian matrices and can also be thought of as
Radon transforms of the shadow. We discuss how the eigenvalues of the sections
are involved in the analysis of the map $\Omega_{N}\rightarrow\dC$
taking $u$ to $(Au,u)$. We remark that, in general, the shadow measure of a nonnormal matrix is absolutely
continuous with respect to area measure on $\dC\equiv\dR^2$ (see [GS2010]).

\subsection{Marginal densities}
Recall that the numerical range of a Hermitian matrix is real and the density
of the shadow is a nonnegative spline function, straightforward to express in
terms of the eigenvalues. This fact can be exploited by means of a type of
Cartesian decomposition.
{\\[.5cm]}
Recall that for $A\in M_{N}$ we define the real part of $A$ by
\[
\mbox{Re}A=\frac{1}{2}\left(  A+A^*\right)  .
\]
Thus $\mbox{Re}A$ is Hermitian, and $A=\mbox{Re}A+i\mbox{Re}\left(  -iA\right)$. We will
be
concerned with the more general $\mbox{Re}\left(  e^{-i
\theta}A\right)  $ with $-\pi\leq\theta\leq\pi$; then $A$ can be expressed as
\beq{E9a.1}{
A=e^{i\theta}\mbox{Re}\left(
e^{-i\theta}A\right)
+ie^{i\theta}\mbox{Re}\left(  -i%
e^{-i\theta}A\right)  .}
\end{equation}
For $-\pi\leq\theta\leq\pi$ let $\ld_k(\theta)$ be the eigenvalues of
Re$(e^{-i\theta}A)$, labeled so that
\[
\ld_1(\theta)\leq\ld_2(\theta)\leq\dots\leq\ld_N(\theta).
\]
Then for each $u\in\Omega_N$ we have $\ld_1(\theta)\leq\mbox{Re}(e^{-i\theta}(Au,u))\leq\ld_N(\theta)$
and
\[
\ld_j(\theta+\pi)=-\ld_{N+1-j}(\theta)\quad(1\leq j\leq N).
\]
A matrix with the property that $\lambda_{1}(  \theta)
<\lambda_{2}(  \theta)  <\dots<\lambda_{N}(
\theta)  $ for $0\leq\theta\leq2\pi$ is called
\textit{generic} in the paper of Jonckheere, Ahmad and Gutkin [JAG1998].
{\\[.5cm]}
We relate the shadow of $\mbox{Re}\left(
e^{-i\theta }A\right)  $ to a marginal density of $P_{A}.$
Write the shadow measure $dP_{A}\left(  z\right)  =p_{A}\left(
z\right)  dm_{2}\left(  z\right)  $ (where $dm_{2}$ is the Lebesgue
measure on $\dC\equiv\dR^2$). Recall: if $f\left(
x,y\right)  $ is a density function on $\dR^2$ with compact
support, then the marginal density along the $x$-axis is
$f_{X}\left(  x\right)  =\int_{-\infty}^{\infty}f\left(  x,y\right)
dy$. Suppose $g\left(  x\right)  $ is a continuous function; then%
\[
E\left[  g\left(  X\right)  \right]  =\int_{-\infty}^{\infty}g\left(
x\right)  f_{X}\left(  x\right)
dx=\int_{-\infty}^{\infty}\int_{-\infty }^{\infty}g\left(  x\right)
f\left(  x,y\right)  dxdy.
\]
Thus the moments, $\int_{-\infty}^{\infty}x^{n}f_{X}\left(  x\right)
dx$ equal $E\left[  X^{n}\right]  $ with respect to the density $f$.
Now replace $x$ by $x\cos\theta+y\sin\theta=\mbox{Re}\left(
e^{-i\theta }z\right)  $ for some fixed $\theta$. The line
orthogonal to $\mbox{Re}\left(
e^{-i\theta}z\right)  =0$ is $\mbox{Re}\left(
ie^{-i\theta}z\right)  =0$. Let
$u=\mbox{Re}\left(  e^{-i\theta}z\right)
,v=\mbox{Re}\left(  ie^{-i\theta}z\right)
$; then $z=e^{i\theta}\left(  u-iv\right)  $ and
so the density of the shadow of $\mbox{Re}\left(
e^{-i\theta}A\right)  $ is
the marginal density of $P_{A}$ for $u$:
\[
f_{U}\left(  u\right)  =\int_{-\infty}^{\infty}p_{A}\left(  e^{i
\theta}\left(  u-iv\right)  \right)  dv.
\]
This is exactly the (2-dimensional) Radon transform of $p_{A}$
evaluated at $\left(  u,\theta\right)
,u\in\dR,-\pi\leq\theta\leq\pi$. We may restate this as follows:
suppose $g\left(  u\right)  $ is real and continuous for
$u\in\dR$; then $E\left[  g\left(  U\right)  \right]  $ with
respect to
$P_{\mbox{Re}\left(  e^{-i\theta}A\right)  }$ is%
\[
\int_{\Omega_{N}}g\left(  \mbox{Re}\left(  e^{-i\theta}%
(Aw,w)\right)  \right)  d\mu(w)     =\int
_{W(A)}g\left(  \frac{1}{2}\left(  e^{-i\theta}z+e^{i%
\theta}\overline{z}\right)  \right)  p_{A}\left(  z\right)
dm_{2}\left(
z\right)
\]
\[
  =\int_{-\infty}^{\infty}g\left(  u\right)  \int_{-\infty}^{\infty}%
p_{A}\left(  e^{i\theta}\left(  u-iv\right)  \right)  dvdu\\
  =E\left[  g\left(  \mbox{Re}\left(
e^{-i\theta}Z\right) \right)  \right]  ,
\]
where the latter expectation is with respect to $P_{A}$.
{\\[.5cm]}
\textbf{Remark:} The Radon transform can be inverted to recover the shadow of $A$
from the shadows of $\mbox{Re}\left(
e^{-i\theta}A\right)$,  $-\pi \leq\theta\leq\pi$. There are
practical algorithms, used in X-ray tomography, which produce
approximations to the inverse transform by using a finite number of
angles $\theta$ (also see (also see [He1984, Ch. 1, Sect. 2]). The Radon transform
approach is worked out thoroughly in [GS2010].
{\\[.5cm]}
For $A\in M_{N}$ let $\xi_{A}\left(  s,t\right)  =\det\left(
I-sA-tA^*\right)  $ and for a Hermitian matrix $H$ let
$\xi_{H}\left(  r\right)
=\det\left(  I-rH\right)  $ (\textquotedblleft$\xi$\textquotedblright%
\ suggests \textquotedblleft characteristic\textquotedblright). For
a power series $h$, $\left[  s^{j}t^{k}\right]  h\left(  s,t\right)
$ denotes the coefficient of $s^{j}t^{k}$ in $h\left(  s,t\right)
$, and $\left[ r^{n}\right]  h\left(  r\right)  $ denotes to
coefficient of $r^{n}$ in $h\left(  r\right)  $,
$j,k,n=0,1,2,\ldots$.

Recall from Proposition 5.3 that the moments of $A$ can be obtained from $\xi_{A}$,%
\[
\nu_{jk}\left(  A\right)  :=\int_{W(A)}z^{j}\overline{z}^{k}%
dP_{A}\left(  z\right)  =\frac{j!k!}{\left(  N\right)  _{j+k}}\left[
s^{j}t^{k}\right]  \xi_{A}\left(  s,t\right)  ^{-1}.
\]
The central moments of a probability distribution are also of
interest. Let $m_{A}=\frac{1}{N}\trace A$, then $E\left[
Z\right]  =m_{A}$. The central
moments can be computed by expanding the integrand in%
\[
\nu_{jk}^{0}\left(  A\right)  :=\int_{W(A)}\left(
z-m_{A}\right) ^{j}\left(  \overline{z}-\overline{m_{A}}\right)
^{k}dP_{A}\left(  z\right) ,
\]
or by using the shifted matrix $A-m_{A}I$.
{\\[.5cm]}
\textbf{Lemma:}
For $A\in M_{N}$ and $c\in\dC$,%
\[
\xi_{A-cI}\left(  s,t\right)  =\left(  1+sc+t\overline{c}\right)
^{N}\xi _{A}\left(
\frac{s}{1+sc+t\overline{c}},\frac{t}{1+sc+t\overline{c}}\right) .
\]
\textbf{Proof:} Indeed,%
\begin{eqnarray*}
\xi_{A-cI}\left(  s,t\right)     &=&\det\left(  I-s\left(
A-cI\right)
-t\left(  A^*-\overline{c}I\right)  \right) \\
&=&\det\left(  \left(  1+sc+t\overline{c}\right)  I-sA-tA^*\right)
  =\left(  1+sc+t\overline{c}\right)  ^{N}\det\left(  I-\frac{sA+tA^*%
}{1+sc+t\overline{c}}\right).
\end{eqnarray*}
QED
{\\[.5cm]}
It is clear that the shadow of $A-cI$ is a translate of $P_{A}$.
Thus
\[
\nu_{jk}^{0}\left(  A\right)  =\frac{j!k!}{\left(  N\right)
_{j+k}}\left[ s^{j}t^{k}\right]  \xi_{A-m_{A}I}\left(  s,t\right)
^{-1},j,k=0,1,2,\ldots.
\]
We consider the (one-dimensional) moments of
$\mbox{Re}\left( e^{-i\theta}A\right)  .$
{\\[.5cm]}
\textbf{Proposition 9.1:}
For
$n=0,1,2,\ldots$,
\[
\int_{-\infty}^{\infty}u^{n}dP_{\mbox{Re}\left(
e^{-i\theta}A\right)  }\left(  u\right)  =\frac{n!}{\left(
N\right) _{n}}\left[  r^{n}\right]  \xi_{A}\left(
\frac{1}{2}re^{-i\theta
},\frac{1}{2}re^{i\theta}\right)  ^{-1}.
\]
\textbf{Proof:} The integral%
\begin{eqnarray*}
\int_{-\infty}^{\infty}u^{n}dP_{\mbox{Re}\left(e^{-i\theta }A\right)  }\left(  u\right)
&=& \frac{1}{2^{n}}\int_{W(A)}\left(e^{-i\theta}z+e^{i\theta}\overline{z}\right)^{n}%
    dP_{A}\left(  z\right) \\
&=& \frac{1}{2^{n}}\sum_{j=0}^{n}\binom{n}{j}e^{i\theta\left(n-2j\right)}
    \int_{W(A)}z^{j}\overline{z}^{n-j}dP_{A}\left(z\right) \\
&=& \frac{1}{2^{n}}\sum_{j=0}^{n}\binom{n}{j}e^{i\theta\left(n-2j\right)}
    \frac{j!\left(  n-j\right)  !}{\left(  N\right)_{n}}\left[s^{j}t^{n-j}\right]  
    \xi_{A}\left(  s,t\right)  ^{-1} \\
&=& \frac{n!}{\left(  N\right)  _{n}}\left[  r^{n}\right] \xi_{A}
    \left(\frac{1}{2}re^{-\mathrm{i}\theta},\frac{1}{2}re^{\mathrm{i}\theta}\right)^{-1}.
\end{eqnarray*}
QED
{\\[.5cm]}
That is, the moments of $P_{\mbox{Re}\left(
e^{-i\theta }A\right)  }$ can be obtained from
$\xi_{\mbox{Re}\left( e^{-i\theta}A\right)  }\left(
r\right)  $. Furthermore, 
$$\xi_{\mbox{Re}\left(
e^{-i\theta}A\right)  }\left(
r\right)  =\det\left(  I-r\mbox{Re}\left(  \left(  e^{-i%
\theta}A\right)  \right)  \right)  =\prod\limits_{j=1}^{N}\left(
1-r\lambda_{j}\left(  \theta\right)  \right).  $$
Here are the basic quantities associated to
$P_{\mbox{Re}\left( e^{-i\theta}A\right)  }$:
{\\[.5cm]}
The mean of $P_{\mbox{Re}\left(  e^{-i\theta}A\right)
}$ is
\[
\frac{1}{2N}\left(  e^{-i\theta}\trace A+e^{i%
\theta}\trace A^*\right)  =\mbox{Re}\left(  e^{-i%
\theta}m_{A}\right)  =\frac{1}{N}\sum_{j=1}^{N}\lambda_{j}\left(
\theta\right).
\]
The variance of $P_{\mbox{Re}\left(  e^{-i\theta
}A\right)  }$ is%
\begin{eqnarray*}
&&\frac{2}{\left(  N\right)  _{2}}\left[  r^{2}\right]  \xi_{\mbox{Re}%
\left(  e^{-i\theta}A\right)  }\left(  r\right)  ^{-1}%
-\mbox{Re}\left(  e^{-i\theta}m_{A}\right)  ^{2}
\\ 
&&=\frac{1}{N\left(  N+1\right)  }\left(
\sum_{j=1}^{N}\lambda_{j}\left( \theta\right)
^{2}-\frac{1}{N}\left(  \sum_{j=1}^{N}\lambda_{j}\left(
\theta\right)  \right)  ^{2}\right)
\\ 
&&=\frac{1}{4N\left(  N+1\right)
}(e^{-2i\theta}\trace\left( \left(  A-m_{A}I\right)
^{2}\right)  +2\trace\left(  \left(
A-m_{A}I\right)  \left(  A^*-\overline{m_{A}}I\right)  \right)
\\ 
&& \phantom{ = } +e^{2i\theta}\trace\left(  \left(  A^*-\overline{m_{A}%
}I\right)  ^{2}\right)  ).
\end{eqnarray*}
Let $\trace\left(  \left(  A-m_{A}I\right)  ^{2}\right)  =ae^{i%
\phi}$ with $a\geq0$; then the variance of
$P_{\mbox{Re}\left( e^{-i\theta}A\right)  }$ is
maximized at $\theta=\frac{\phi}{2}$ and minimized at
$\theta=\frac{\phi}{2}\pm\frac{\pi}{2}$. There is a relation with
the 2-dimensional variance of $P_{A}$, namely,
\[
\int_{W(A)}\left\vert z-m_{A}\right\vert ^{2}dP_{A}\left(
z\right)  =\frac{1}{N\left(  N+1\right)
}\trace\left(  \left(  A-m_{A}I\right)  \left(  A^*-\overline{m_{A}%
}I\right)  \right).
\]
The central moments of $\mbox{Re}\left(e^{-i\theta }A\right) $ 
can be obtained from
\begin{eqnarray*}
\xi_{\mbox{Re}\left( e^{-i\theta}A\right)
-\mbox{Re}\left(  e^{-i\theta }m_{A}\right)
I}\left(  r\right)
  =\left(  1+\mbox{Re}\left(
e^{-i\theta}m_{A}\right)  r\right)
^{N}\xi_{\mbox{Re}\left( e^{-i\theta}A\right)
}\left(  \dfrac{r}{1+r\mbox{Re}\left(
e^{-i\theta}m_{A}\right)  }\right).
\end{eqnarray*}

\subsection{The shadow of a Hermitian matrix}

The density
function for $P_{H}$ is simple to find, given the eigenvalues of a
Hermitian matrix $H$. Suppose $H$ is not scalar; then $H$ has at
least two different eigenvalues and the shadow is absolutely
continuous on $\dR$. Let $dP_{H}\left( x\right) =p_{H}\left(
x\right)  dx$. The following is the basic fact.
{\\[.5cm]}
\textbf{Lemma 9.2:} Suppose $1\leq m<N$ and $\lambda\neq0$. For $n=0,1,2,\ldots$%
\[
m\binom{N-1}{m}\lambda^{1-N}\int_{0}^{\lambda}x^{n}x^{m-1}\left(
\lambda-x\right)  ^{N-m-1}dx    =\frac{\left(  m\right)
_{n}}{\left(
N\right)  _{n}}\lambda^{n}\quad(\lambda>0);
\]\[
m\binom{N-1}{m}\left(  -\lambda\right)
^{1-N}\int_{\lambda}^{0}x^{n}\left( -x\right)  ^{m-1}\left(
x-\lambda\right)  ^{N-m-1}dx    =\frac{\left( m\right)
_{n}}{\left(  N\right)  _{n}}\lambda^{n}\quad(\lambda<0).
\]
To express the density functions for all real arguments we use the notation%
\[
x_{+}=\max\left(  x,0\right)  ,
\]
with the convention that $x_{+}^{0}=1$ for $x\geq0$ and $=0$ for
$x<0$. Thus
the density in the first part of the lemma equals $m\binom{N-1}{m}%
\lambda^{1-N}x_{+}^{m-1}\left(  \lambda-x\right)  _{+}^{N-m-1}$ for
$x\in\dR$.
{\\[.5cm]}
Suppose $H$ is Hermitian, not a multiple of $I$, and $\xi_{H}\left(
r\right)
=\prod\limits_{i=1}^{N}\left(  1-r\lambda_{i}\right)  =\prod\limits_{i=1}%
^{M}\left(  1-r\mu_{i}\right)  ^{m_{i}}$, where $\left\{
\mu_{1},\ldots ,\mu_{M}\right\}  $ is the set of distinct nonzero
eigenvalues of $H$ and $\sum_{i=1}^{M}m_{i}\leq N$. By hypothesis,
$H$ has at least two different eigenvalues (each $m_{i}<N$). There
are are unique real numbers $\beta
_{ij},1\leq i\leq M,1\leq j\leq m_{i}$ such that%
\[
\frac{1}{\xi_{H}\left(  r\right)  }=\sum_{i=1}^{M}\sum_{j=1}^{m_{i}}%
\frac{\beta_{ij}}{\left(  1-r\mu_{i}\right)
^{j}}=\sum_{i=1}^{M}\sum
_{j=1}^{m_{i}}\beta_{ij}\sum_{n=0}^{\infty}\frac{\left(  j\right)  _{n}}%
{n!}\mu_{i}^{n}r^{n},\left\vert r\right\vert
<\frac{1}{\min_{i}\left\vert \mu_{i}\right\vert }.
\]
Thus for $n=0,1,2,\ldots$%
\[
\int_{-\infty}^{\infty}x^{n}dP_{H}\left(  x\right)
=\sum_{i=1}^{M}\sum
_{j=1}^{m_{i}}\beta_{ij}\frac{\left(  j\right)  _{n}}{\left(  N\right)  _{n}%
}\mu_{i}^{n}.
\]
By the lemma, the density function of $P_{H}$ is%
\begin{eqnarray*}
p_{H}\left(  x\right)   &=& \sum_{\mu_{i}>0}\sum_{j=1}^{m_{i}}\mu_{i}%
^{1-N}\beta_{ij}j\binom{N-1}{j}x_{+}^{j-1}\left(  \mu_{i}-x\right)
_{+}^{N-j-1}\quad \\
& & +  \sum_{\mu_{i}<0}\sum_{j=1}^{m_{i}}\left(  -\mu_{i}\right)  ^{1-N}%
\beta_{ij}j\binom{N-1}{j}\left(  -x\right)  _{+}^{j-1}\left(  x-\mu
_{i}\right)  _{+}^{N-j-1}.
\end{eqnarray*}

\subsection{The critical curves in W(A)}

Suppose that in some interval $\theta_{1}<\theta<\theta_{2}$ the
eigenvalues of $\mbox{Re}\left( e^{-i\theta}A\right)  $ are pairwise
distinct, and there are eigenvectors $\psi^{\left( j\right)  }\left(
\theta\right)  $ so that%
\[
\mbox{Re}\left(  e^{-i\theta}A\right)  \psi^{\left( j\right)
}\left(  \theta\right)     =\lambda_{j}\left( \theta\right)
\psi^{\left(  j\right)  }\left(  \theta\right)  ,
\]
where
\[
\left\vert \psi^{\left(  j\right)  }\left(  \theta\right)
\right\vert   =1\quad (1\leq j\leq N).
\]
The image $\left\{  (A\psi^{\left(  j\right)  }\left(  \theta\right),
\psi^{\left(  j\right)  }\left(  \theta\right))
:\theta_{1}<\theta <\theta_{2}\right\}  $ is called a critical curve
(see [JAG1998,Theorem 5, p. 238]).
{\\[.5cm]}
\textbf{Lemma 9.3:} Suppose $H\left(  \theta\right)  ,\psi\left(  \theta\right)
,\lambda\left( \theta\right)  $ are differentiable functions on
$\theta_{1}<\theta<\theta _{2}$ such that $H\left(  \theta\right)  $
is Hermitian, $\psi\left( \theta\right)  \in\Omega_{N}$,
$\lambda\left(  \theta\right) \in\dR$, and $H\left(
\theta\right)  \psi\left(  \theta\right) =\lambda\left(
\theta\right)  \psi\left(  \theta\right)  $; then $\frac
{d}{d\theta}\lambda\left(  \theta\right)  =((  \frac{d}{d\theta}H\left(  \theta\right))
\psi\left( \theta\right),\psi(\theta))  $.
{\\[.5cm]}
\textbf{Proof:} Write $\lambda\left(  \theta\right)
=\psi\left( \theta\right)  ^*H\left(  \theta\right)  \psi\left(
\theta\right)  $ and differentiate to obtain%
\begin{eqnarray*}
\frac{d}{d\theta}\lambda\left(  \theta\right) &=& 
  \psi\left(\theta\right) ^*\left(  \frac{d}{d\theta}H\left(  \theta\right)
  \right)  \psi\left( \theta\right)  +\frac{d}{d\theta}\left(
  \psi\left(  \theta\right)  ^*\right)  H\left(  \theta\right)
  \psi\left(  \theta\right)  +\psi\left( \theta\right)  ^*H\left(
  \theta\right)  \frac{d}{d\theta}\psi\left(
  \theta\right)
\\ 
&=&\psi\left(  \theta\right)  ^*\left(
  \frac{d}{d\theta}H\left( \theta\right)  \right)  \psi\left(
  \theta\right)  +\lambda\left( \theta\right)  \left(
  \frac{d}{d\theta}\left(  \psi\left(  \theta\right) ^*\right)
  \psi\left(  \theta\right)  +\psi\left(  \theta\right)  ^*
  \frac{d}{d\theta}\psi\left(  \theta\right)  \right)
\\
&=&\psi\left(  \theta\right)  ^*\left(
  \frac{d}{d\theta}H\left( \theta\right)  \right)  \psi\left(
  \theta\right)  ,
\end{eqnarray*}
because $\psi\left(  \theta\right)  ^*\psi\left(  \theta\right)
=1.$
QED
{\\[.5cm]}
\textbf{Proposition 9.4:} For $\theta_{1}<\theta<\theta_{2}$ the critical curve satisfies\\
$\bigl(A\psi^{(j)}(\theta),\psi^{(j)}(\theta)\bigr)=$ $\psi^{\left( j\right)  }\left(  \theta\right)  ^*A\psi^{\left(
j\right)  }\left( \theta\right)  =e^{i\theta}\left(
\lambda_{j}\left(  \theta\right) +i\lambda_{j}^{\prime}\left(
\theta\right)  \right)  $.
{\\[.5cm]}
\textbf{Proof:} Let $H\left(  \theta\right)  =\mbox{Re}\left( e^{-i\theta
}A\right)  =\left(  \cos\theta\right)  A_{1}+\left(  \sin\theta\right)  A_{2}%
$, where $A_{1}=\frac{1}{2}\left(  A+A^*\right)  $ and
$A_{2}=\frac
{1}{2i}\left(  A-A^*\right)  $. Thus%
\[
\psi^{\left(  j\right)  }\left(  \theta\right)  ^*A\psi^{\left(
j\right)  }\left(  \theta\right)  =\psi^{\left(  j\right)  }\left(
\theta\right)  ^*A_{1}\psi^{\left(  j\right)  }\left(
\theta\right)
+i\psi^{\left(  j\right)  }\left(  \theta\right)  ^*A_{2}%
\psi^{\left(  j\right)  }\left(  \theta\right)  .
\]
By definition and the lemma%
\[
\psi^{\left(  j\right)  }\left(  \theta\right)  ^*\bigl(
\left( \cos\theta\right)  A_{1}+\left(  \sin\theta\right)
A_{2}\bigr) \psi^{\left(  j\right)  }\left(  \theta\right)
=\lambda_{j}\left(
\theta\right)  ,
\]
and
\[
\psi^{\left(  j\right)  }\left(  \theta\right)  ^*\bigl(
-\left( \sin\theta\right)  A_{1}+\left(  \cos\theta\right)
A_{2}\bigr) \psi^{\left(  j\right)  }\left(  \theta\right)
=\lambda_{j}^{\prime }\left(  \theta\right)  .
\]
These equations are easily solved to establish the formula.
QED
{\\[.5cm]}
The analyticity of $\lambda_{j}\left(  \theta\right)  $ is shown in
[JAG1998, Lemma 2, p. 240].
{\\[.5cm]}
To get an idea of the structure of critical curves one has to
distinguish the generic and non-generic cases. If $A$ is generic
then $\lambda_{j}\left( \theta+\pi\right)  =-\lambda_{N+1-j}\left(
\theta\right)  $ for $1\leq j\leq N$, the curve $C_{j}:=\left\{
e^{i\theta}\left(  \lambda_{j}\left( \theta\right)
+i\lambda_{j}^{\prime}\left(  \theta\right)  \right)
:0\leq\theta\leq2\pi\right\}  $ agrees with $C_{N+1-j}$ (as a
point-set) so there are $\left\lfloor \frac{N+1}{2}\right\rfloor $
critical curves (see [JAG1998, Theorem 13, p. 244]); the outside $C_{1}$
is the boundary of $\Lambda_{A}$. Example 9.5.3 below is a
generic $3\times3$ matrix.
{\\[.5cm]}
When using numerical techniques for solving the characteristic
equation for some number of angles (for example
$\theta=j\pi/m,j=0,\ldots,m$) the value of $\lambda^{\prime}\left(
\theta\right)  $ can be computed as follows: let
$p\left(  \theta,\lambda\right)  =\det\left(  \lambda I-\mbox{Re}%
\left(  e^{-i\theta}A\right)  \right)  $, differentiate the equation
$p\left(  \theta,\lambda\left(  \theta\right)  \right)  =0$ to
obtain $$\frac{\partial}{\partial\theta}p\left( \theta,\lambda\left(
\theta\right)
\right)  +\lambda^{\prime}\left(  \theta\right)  \frac{\partial}%
{\partial\lambda}p\left(  \theta,\lambda\left(  \theta\right)
\right)  =0$$ (so the value of $\lambda\left(  \theta\right)  $
determines $\lambda^{\prime }\left(  \theta\right)  $ except
possibly for isolated points where
$\frac{\partial}{\partial\lambda}p\left(  \theta,\lambda\left(
\theta\right) \right)  =0$; this indicates repeated roots which do
not occur in the generic case).
{\\[.5cm]}
In the non-generic case the same critical curve can arise from
different eigenvalues: let $\theta_{0}$ be an angle for which the
eigenvalues are all distinct and ordered by $\lambda_{1}\left(
\theta_{0}\right)  <\ldots <\lambda_{N}\left(  \theta_{0}\right)  $;
consider each $\lambda_{i}\left( \theta\right)  $ as a real-analytic
function in $\theta$ and extend it to the interval
$\theta_{0}\leq\theta\leq\theta_{0}+\pi$. Because this is the
non-generic case the curves $\lambda_{i}\left(  \theta\right)  $ may
cross in the open interval (a finite number of times by
analyticity). Form a set-partition of $\left\{  1,2,\ldots,N\right\}
$ by declaring $i$ and $j$ equivalent if $\lambda_{i}\left(
\theta_{0}\right)  =-\lambda_{j}\left( \theta_{0}+\pi\right)  $; the
relation is extended by transitivity. The equivalence classes
correspond to distinct critical curves. There may be only one class;
consider Example 9.5.2. In this case the boundary of
$W(A)$ is the convex hull of the outside critical curve (from
the class containing $1$).
{\\[.5cm]}
In the situation of radially symmetric shadows (see section 8) the
critical curves are circles centered at the origin.

\subsection{A geometric approach}

Any matrix can be expressed as a sum of two matrices with orthogonal
1-dimensional numerical ranges. For a fixed $\theta$ with $-\frac{\pi}{2}%
\leq\theta\leq\frac{\pi}{2}$ we can write%
\[
A=e^{i\theta}\operatorname{Re}\left(  e^{-i\theta}A\right)
+ie^{i\theta}\operatorname{Re}\left(  e^{-i\left(
\theta+\pi/2\right)  }A\right)
\]
Let $U_{1},U_{2}\in\mathcal{U}\left(  N\right)  $ 
 satisfy $U_{1}^{*}\operatorname{Re}\left(  e^{-i\theta}A\right)  U_{1}=B_{1}%
,U_{2}^*\operatorname{Re}\left(  e^{-i\left(  \theta
+\pi/2\right)  }A\right)  U_{2}=B_{2}$ where $B_{1},B_{2}$ are diagonal
matrices, such that $\left(  B_{1}\right)  _{jj}=\lambda_{j}\left(
\theta\right)  $ and $\left(  B_{2}\right)  _{jj}=\lambda_{j}\left(
\theta+\frac{\pi}{2}\right)  $ for $1\leq j\leq N$. Then for $\psi\in
\Omega_{N}$%
\[
\psi^*A\psi=e^{i\theta}\psi^*A_{1}\psi+i%
e^{i\theta}\psi^*A_{2}\psi=e^{i\theta}\psi^*%
U_{1}B_{1}U_{1}^*\psi+ie^{i\theta}\psi^*U_{2}%
B_{2}U_{2}^*\psi.
\]
By the unitary invariance of the range (and the shadow) we may replace
(generic) $\psi$ by (generic) $U_{1}\psi$. Thus%
\begin{equation}
\left(  U_{1}\psi\right)  ^*A\left(  U_{1}\psi\right)  =e^{i%
\theta}\sum_{j=1}^{N}\lambda_{j}\left(  \theta\right)  \left\vert \psi
_{j}\right\vert ^{2}+ie^{i\theta}\sum_{j=1}^{N}\lambda
_{j}\left(  \theta+\frac{\pi}{2}\right)  \left\vert \left(  U_{2}^*%
U_{1}\psi\right)  _{j}\right\vert ^{2}. \label{decomp}%
\end{equation}
The value remains unchanged if $U_{2}^*U_{1}$ is replaced by
$U=D_{2}U_{2}^*U_{1}D_{1}$ where $D_{1},D_{2}$ are arbitrary diagonal
unitary matrices (in $M_N$). For example, choose
$D_{1},D_{2}$ so that $U_{1,j}\geq0$ and $U_{j,1}\geq0$ for $1\leq j\leq N$.
Thus the numerical range and shadow can be interpreted in terms of a mapping
from
\[
T_{N-1}\left[  U\right]  :=\left\{  \left(  \left(  \left\vert \psi
_{j}\right\vert ^{2}\right)  _{j=1}^{N},\left(  \left\vert U\psi
_{j}\right\vert ^{2}\right)  _{j=1}^{N}\right)  :\psi\in\Omega_{N}\right\}
\subset \Delta_N\times \Delta_N%
\]
to $\dC$. Every vector $\left(  \left\vert \psi_{j}\right\vert
^{2}\right)  _{j=1}^{N}$ appears as a first and as a second component of a
point in $T_{N-1}\left[  U\right]  $. The unitarily invariant measure on
$\Omega_{N}$ induces a measure on $T_{N-1}\left[  U\right]  $, and the shadow
of $A$ is the image of this measure under the map%
\[
\left(  \left(  \left\vert \psi_{j}\right\vert ^{2}\right)  _{j=1}^{N},\left(
\left\vert U\psi_{j}\right\vert ^{2}\right)  _{j=1}^{N}\right)  \mapsto
e^{i\theta}\sum_{j=1}^{N}\left(  \lambda_{j}\left(  \theta\right)
\left\vert \psi_{j}\right\vert ^{2}+i\lambda_{j}\left(  \theta
+\frac{\pi}{2}\right)  \left\vert \left(  U\psi\right)  _{j}\right\vert
^{2}\right)  .
\]
The case $N=2$ can be explicitly described. Let
\[
U=%
\begin{pmatrix}
\cos\theta_{0} & \sin\theta_{0}\\
\sin\theta_{0} & -\cos\theta_{0}%
\end{pmatrix}
,\psi=%
\begin{bmatrix}
e^{i\phi_{1}}\cos\theta_{1}\\
e^{i\phi_{2}}\sin\theta_{1}%
\end{bmatrix}
,
\]
with $0\leq\theta_{0},\theta_{1}\leq\frac{\pi}{2}$ and $-\pi\leq\phi_{1}%
,\phi_{2}\leq\pi$. It suffices to consider $U$ of this form ($\theta_{0}$ is
fixed). Then%
\[
\left\vert U\psi_{1}\right\vert ^{2}    =\frac{1}{2}+\frac{1}{2}\cos
2\theta_{1}\cos2\theta_{0}+\frac{1}{2}\sin2\theta_{1}\sin2\theta_{0}%
\cos\left(  \phi_{1}-\phi_{2}\right)  ,
\]
and
\[
\left\vert U\psi_{2}\right\vert ^{2}    =\frac{1}{2}-\frac{1}{2}\cos
2\theta_{1}\cos2\theta_{0}-\frac{1}{2}\sin2\theta_{1}\sin2\theta_{0}%
\cos\left(  \phi_{1}-\phi_{2}\right)  .
\]
Thus $T_{1}\left[  U\right]$ is
\[
\Bigl\{\left(  \frac{1}{2}+\frac{1}{2}x_{1}\left(
\theta_{1}\right)  ,\frac{1}{2}-\frac{1}{2}x_{1}\left(  \theta_{1}\right)
\right)  ,\left(  \frac{1}{2}+\frac{1}{2}x_{2}\left(  \theta_{1},\theta
_{2},\phi\right)  ,\frac{1}{2}-\frac{1}{2}x_{2}\left(  \theta_{1},\theta
_{2},\phi\right)  \right)  :
\]\[
0    \leq\theta_{1}\leq\frac{\pi}{2},0\leq\phi\leq\pi \Bigr\},
\]
where
$x_{1}\left(  \theta_{1}\right)     =\cos2\theta_{1}$ and
\[
x_{2}\left(  \theta_{1},\theta_{0},\phi\right)     =\cos2\theta_{1}%
\cos2\theta_{0}+\sin2\theta_{1}\sin2\theta_{0}\cos\phi.
\]
As expected (compare section 2) this forms an ellipse (including the interior). Changing
coordinates we transform  to the square $\left\{  \left(
x_{1},x_{2}\right)  :-1\leq x_{1},x_{2}\leq1\right\}  $; then $T_{1}\left[
U\right]  $ maps to $\left\{  \left(  x_{1},x_{2}\right)  :x_{1}%
^{2}-2x_{1}x_{2}\cos2\theta_{0}+x_{2}^{2}\leq\sin^{2}2\theta_{0}\right\}  $.
In the degenerate normal case this reduces to the interval $\left\{  \left(
x_{1},x_{1}\right)  :-1\leq x_{1}\leq1\right\}  $.
{\\[.5cm]}
The invariant measure on $\Omega_{2}$ is $\frac{1}{2\pi^{2}}\sin\theta_{1}%
\cos\theta_{1}d\theta_{1}d\phi_{1}d\phi_{2}$. This is mapped to the measure
$$
\frac{1}{2\pi}\left(  \sin^{2}2\theta_{0}-x_{1}^{2}+2x_{1}x_{2}\cos
2\theta_{0}-x_{2}^{2}\right)  ^{-\frac{1}{2}}dx_{1}dx_{2}.
$$

\subsection{Examples}

Example 9.5.1 Let
\[
A_1=%
\begin{pmatrix}
0 & 1 & 0 & 0\\
0 & 0 & 1 & 0\\
0 & 0 & 0 & 1\\
0 & 0 & 0 & 0
\end{pmatrix}
.
\]
Then $\xi_{A_1}\left(  s,t\right)  =1-3st+s^{2}t^{2}$ and $\xi
_{\mbox{Re}\left(  e^{-i\theta}A_1\right)  }\left(  r\right)
=1-\frac{3}{4}r^{2}+\frac{1}{16}r^{4}$. The eigenvalues of $\mbox{Re}%
\left(  e^{-i\theta}A_1\right)  $ are $\frac{1}{4}\left(  \pm1\pm
\sqrt{5}\right)  $, independent of $\theta$, labeled so that $\lambda
_{1}<\lambda_{2}<0<\lambda_{3}<\lambda_{4}$.  Thus the density for
Re$(e^{-i\theta}A_1)$ is given by
\[
P_1(x)=
6\left(  1-\frac{\sqrt{5}}{5}\right)  \times\left(  \left(  -x\right)
_{+}^{0}\left(  x-\lambda_{1}\right)  _{+}^{2}+x_{+}^{0}\left(  \lambda
_{4}-x\right)  _{+}^{2}\right)\,\,\,-
\]\[
6\left(  1+\frac{\sqrt{5}}{5}\right)  \left(  \left(  -x\right)  _{+}%
^{0}\left(  x-\lambda_{2}\right)  _{+}^{2}+x_{+}^{0}\left(  \lambda
_{3}-x\right)  _{+}^{2}\right)  .
\]
In fact the shadow has circular symmetry, as discussed in sections 7 and 8. The
critical curves are $z=\frac{1}{4}\left(  1+\sqrt{5}\right)  e^{i%
\theta}$ and $z=\frac{1}{4}\left(  \sqrt{5}-1\right)  e^{i\theta
},0\leq\theta\leq2\pi$.
{\\[.5cm]}
Example 9.5.2
Let%
\[
A_2=%
\begin{pmatrix}
0 & 1 & 1 & 1\\
0 & 0 & 1 & 1\\
0 & 0 & 0 & 1\\
0 & 0 & 0 & 0
\end{pmatrix}
.
\]
Then
\[
\xi_{A_2}\left(  s,t\right)     =1-6st-4st\left(  s+t\right)  -st\left(
s^{2}+st+t^{2}\right)
\]
and
\[
\xi_{\mbox{Re}\left(  e^{-i\theta}A_2\right)  }\left(
r\right)     =\left(  1+r-r^{2}\left(  \frac{1}{4}-\frac{1}{2}\cos
\theta\right)  \right)  \left(  1-r-r^{2}\left(  \frac{1}{4}+\frac{1}{2}%
\cos\theta\right)  \right).
\]
\begin{figure}[h!]
\includegraphics[width=\figurewidth]{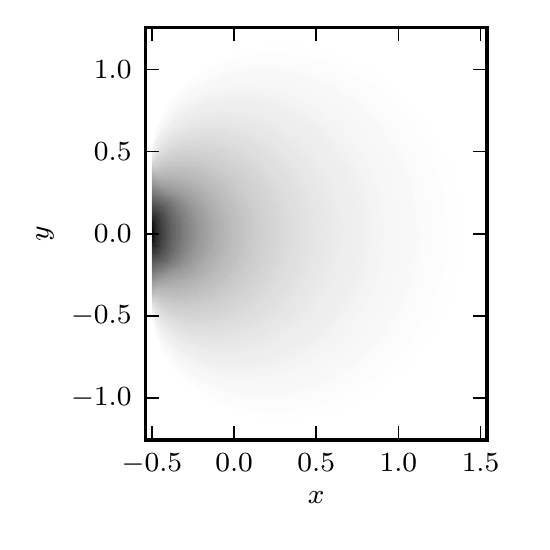} 
\includegraphics[width=\figurewidth]{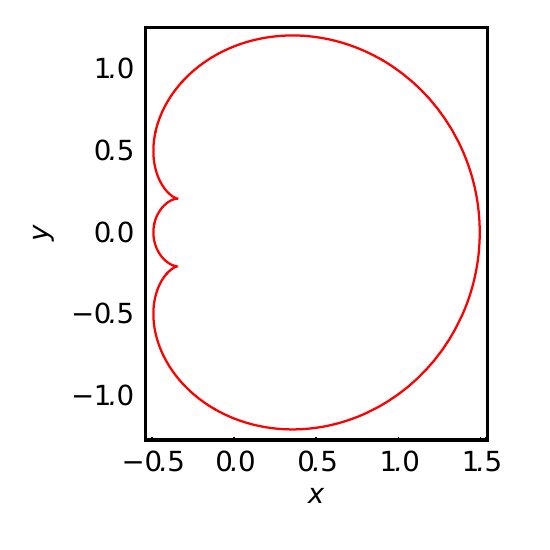}
\caption{Numerical shadow for matrix $A_2$ of size $N=4$ 
             and the corresponding critical lines.}
\label{fig:ex2a}
\end{figure}
The eigenvalues of $\mbox{Re}\left(  e^{-i\theta}A_2\right)  $
are $\frac{1}{2}\pm\cos\frac{\theta}{2},-\frac{1}{2}\pm\sin\frac{\theta}{2}$.
The eigenvalues are $\left[  \frac{3}{2},-\frac{1}{2},-\frac{1}{2},-\frac
{1}{2}\right]  $ at $\theta=0$ and $\pm\frac{1}{2}\pm\frac{\sqrt{2}}{2}$ at
$\theta=\frac{\pi}{2}$. In the range $-\frac{\pi}{2}\leq\theta\leq\frac{\pi
}{2}$ we have%
\[
\lambda_{4}\left(  \theta\right)     =\frac{1}{2}+\cos\frac{\theta}{2},\quad
\lambda_{1}\left(  \theta\right)     =-\frac{1}{2}-\left\vert \sin
\frac{\theta}{2}\right\vert,
\]
and so $-\frac{1}{2}-\left\vert \sin\frac{\theta}{2}\right\vert \leq
\mbox{Re}\left(  e^{-i\theta}\psi^*A_2\psi\right)
\leq\frac{1}{2}+\cos\frac{\theta}{2}$ for $\psi\in\Omega_{N}$. The triple
eigenvalue $-\frac{1}{2}$ at $\theta=0$ results in a pronounced peak in
the density for
$\mbox{Re}(A_2)$:%
\[
P_2(x)=
6\left(
\frac{1}{2}+x\right)  _{+}\left(  -x\right)  _{+}^{2}+9\left(  \frac{1}%
{2}+x\right)  _{+}\left(  -x\right)  _{+}+\frac{27}{8}\left(  \frac{1}%
{2}+x\right)  _{+}^{2}\left(  -x\right)  _{+}^{0}
\]\[
  +\frac{3}{8}x_{+}^{0}\left(  \frac{3}{2}-x\right)  _{+}^{2}.
\]
The matrix $A_2$ is non-generic and there is only one critical curve
which  has two cusps as shown in Fig. \ref{fig:ex2a}.
In this example the boundary of the shadow is the convex hull of the critical curve, 
so the line segment $[(-1-i)/2, (-1+i)/2]$ is a part of the boundary.
One representation of critical lines is $z=e^{i\theta}\left(  \frac{1}{2}%
+\cos\frac{\theta}{2}-\frac{i}{2}\sin\frac{\theta}{2}\right)
,0\leq\theta\leq4\pi$. The line segment joining $\frac{-1-i}{2}$ to
$\frac{-1+i}{2}$ is part of the boundary of $\Lambda_{A_2}$. 
The cusps are at $\dfrac{-19\pm5\sqrt{5}}{54}$.
{\\[.5cm]}
Example 9.5.3
Let
\[
A_3=%
\begin{pmatrix}
0 & 1 & 1\\
0 & i & 1\\
0 & 0 & -1
\end{pmatrix}
.
\]

Then $m_{A_3}=\frac{1}{3}\left(  -1+i\right)  $,
\begin{eqnarray*}
\xi_{A_3}\left(  s,t\right)
&=& 1+\left(  1-i\right)  s+\left(1+i\right)  t-i
  \left(  s^{2}-t^{2}\right)  -3st-st\left(
  \left(  2-i\right)  s+\left(  2+i\right)  t\right), \\
\xi_{A_3-m_{A_3}I}\left(  s,t\right)
&=& 1-\frac{13}{3}st-\frac{i}{3}\left(  s^{2}-t^{2}\right) 
  +\frac{5}{27}\left(  \left(  1+i\right)  s^{3}+\left(  1-i\right)  
  t^{3}\right)\\
&& -\frac{st}{9}\left(\left(  8+i\right)  s+\left(  8-i\right)  t\right),
\end{eqnarray*}
and
\[
\xi_{\mbox{Re}\left(  e^{-i\theta}A_3\right)  }\left(
r\right)     =1+\left(  \cos\theta-\sin\theta\right)  r-\left(  \frac{3}%
{4}+\sin\theta\cos\theta\right)  r^{2}+\left(  \frac{1}{4}\sin\theta-\frac
{1}{2}\cos\theta\right)  r^{3}.
\]
\begin{figure}[ht!]
\includegraphics[width=\figurewidth]{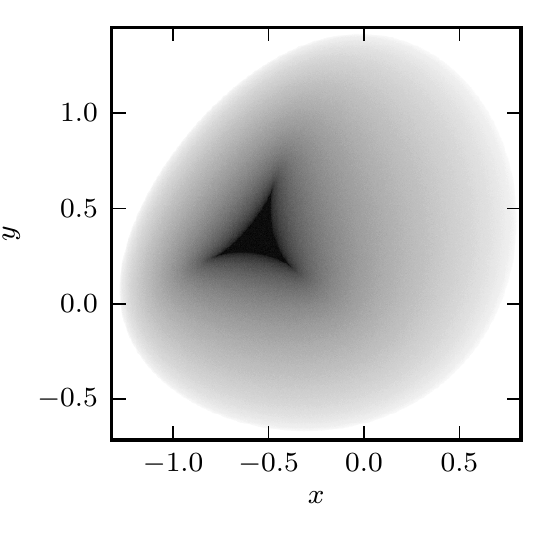} 
\includegraphics[width=\figurewidth]{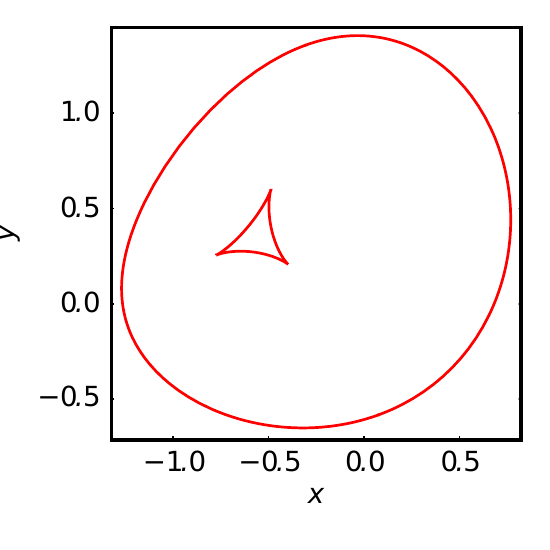}
\caption{Numerical shadow for matrix $A_3$ of size $N=3$ and the corresponding critical lines.}
\label{fig:ex2}
\end{figure}
The variance of $P_{\mbox{Re}\left(  e^{-i\theta}A_3\right)  }$
equals $\frac{1}{72}\left(  13+2\sin2\theta\right)  $. The eigenvalues of
$\mbox{Re}\left(  e^{-i\theta}A_3\right)  $ can be found
approximately, or analytically by the classical formula. For $\theta=0$ the
eigenvalues are $\left[  -\frac{1+\sqrt{17}}{4},-\frac{1}{2},\frac
{-1+\sqrt{17}}{4}\right]  $ for $\theta=0$ and $\left[
-0.6715,0.2647,1.407\right]  $ (rounded) for $\theta=\frac{\pi}{2}$. The
matrix is generic; there are two critical curves: one is the boundary of
$\Lambda_{A_3}$ and the other is a triangular curve with three cusps.
A comparison of the numerical shadow for this matrix and its critical
lines is presented in Fig. \ref{fig:ex2}.
Since $N=3$ the density for $\mbox{Re}\left(  e^{-i\theta
}A_3\right)  $ is piecewise linear:
\[
P_3(x)=
\left(
-x\right)  _{+}^{0}\left(  x+\frac{1+\sqrt{17}}{4}\right)  _{+}-2\left(
-x\right)  _{+}^{0}\left(  x+\frac{1}{2}\right)  _{+}
\]\[
  +\left(  1-\frac{1}{\sqrt{17}}\right)  x_{+}^{0}\left(  \frac{\sqrt{17}%
-1}{4}-x\right)  _{+}.
\]

\section{Direct sums (block diagonal matrices)}

\begin{figure}[ht!]
  \begin{center}
    \includegraphics[width=\figurewidth]{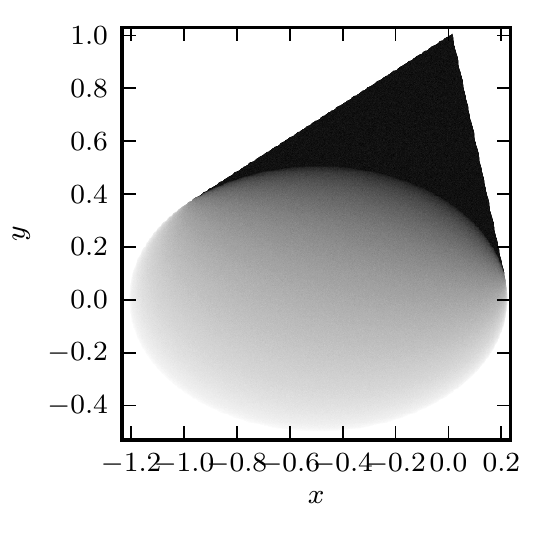}
    \caption[]{
    Shadow of block diagonal matrix
    $
    \left(
    \begin{smallmatrix}
     -1 & 0 \\
     1 & 0
    \end{smallmatrix}
    \right)
    \oplus
    \left(
    \begin{smallmatrix}
    i
    \end{smallmatrix}
    \right)
    $
    .
    }
    \label{sec2fig1}
  \end{center}
\end{figure}

Concerning the direct sum $A\oplus B$
(block diagonal matrix) of matrices $A,B$, it is well--known that
\[
W(A\oplus B)=\conv\{W(A)\cup W(B)\}
\]
(see, for example, [B1997, Exercise I.3.1]). In our context it is natural to ask how the
numerical shadow of $A\oplus B$ is distributed over $\conv\{W(A)\cup W(B)\}$. Here $A$
and $B$ may be of different sizes - see an example presented in Fig. \ref{sec2fig1}.
We consider then $A\oplus B\in M_N$ with $A\in M_n$ and
$B\in M_m$, so that $n+m=N$. Given $u\in\Omega_N$ (distributed according to $\mu$, as usual),
let $u=v_1\oplus v_2$ where $v_1\in\dC^n$ and $v_2\in\dC^m$; then $\|v_1\|^2+\|v_2\|^2=1$.
It is known that $t=\|v_1\|^2$ has a beta--density given by
\beq{E10.1}{
q(t)=\frac{(n+m-1)!}{(n-1)!\; (m-1)!} t^{n-1}(1-t)^{m-1}\quad (t\in[0,1]).}
\end{equation}
From this one can deduce that the shadow measure $P_{A\oplus B}$ is an
``$(n,m)$--beta mixture'' of the shadow measures $P_A$ and $P_B$. Compare
[GS2010, section 2.2].

Another version of this result relates the densities corresponding to $A$, $B$, and $A\oplus B$:
{\\[.5cm]}
\textbf{Proposition 10.1:} If $p_A(z),p_B(z)$ are the shadow densities for $A,B$, then the corresponding
density $p$ for $A\oplus B$ is given by
\beq{E10.2}{
p(z)=\int_0^1 q(t)\Bigl(\int_\dC t^{-2}p_A((z-w)/t)(1-t)^{-2}p_B\bigl(w/(1-t)\bigr)\,dw\Bigr)\,dt,}
\end{equation}
where $q(t)$ is as in (\ref{E10.1}).\\
\textbf{Proof:} For $u\in\Omega_N$ we have
\[
\bigl((A\oplus B)u,u\bigr)=(Av_1,v_1)+(Bv_2,v_2)
\]\[
=\|v_1\|^2(Au_1,u_1)+\|v_2\|^2(Bu_2,u_2),
\]
where $u_j=v_j/\|v_j\|$. Note that $u_1\in\Omega_n$, $u_2\in\Omega_m$, $u_1,u_2$ are stochastically independent,
and they have the corresponding uniform distributions over $\Omega_n.\Omega_m$.
{\\[.5cm]}
Hence $((A\oplus B)u,u)=tZ_1+(1-t)Z_2$ where $Z_1$ and $Z_2$ are independent complex random variables
with densities $p_A(z)$ and $p_B(z)$. Thus
\[
p(z)=\int_0^1 q(t)g(z,t)\,dt
\]
where $g(z,t)$ is the density of the independent (for each fixed $t$) sum $tZ_1+(1-t)Z_2$.
This density is given by the usual convolution formula
\[
g(z,t)=\int_\dC g_1(z-w)g_2(w)\,dw,
\]
where $g_1$ is the density of $tZ_1$ and $g_2$ is the density of $(1-t)Z_2$. If a complex random
variable $Z$ has density $h(z)$ with respect to area on $\dC$, then $tZ$ (where $t\in\dR$) has density
$t^{-2}h(z/t)$. Hence $g_1(z-w)=t^{-2}p_A((z-w)/t)$, $g_2(w)=(1-t)^{-2}p_B(w/(1-t))$ and (\ref{E10.2})
follows. QED

\section{Zernike expansions}

Given our
methods for evaluating the moments of shadow measures (see section 5),
it is natural to construct orthogonal polynomial approximations using Zernike
polynomials. These provide one way to generate pictures of specific
numerical shadows.
{\\[.5cm]}
The complex Zernike polynomials $Z_{mn}\left(  z,\overline{z}\right)  $ are
orthogonal for area measure on the unit disk $\left\{  z\in\mathbb{C}%
:\left\vert z\right\vert \leq1\right\}  $. They can be defined by%
\[
Z_{mn}\left(  z,\overline{z}\right)     =z^{m-n}\sum_{j=0}^{n}\frac{\left(
m+n-j\right)  !}{\left(  m-j\right)  !\left(  n-j\right)  !j!}\left(
-1\right)  ^{j}\left(  z\overline{z}\right)  ^{n-j}\quad (m\geq n)
\]\[
Z_{mn}\left(  z,\overline{z}\right)     =\overline{z}^{n-m}\sum_{j=0}%
^{m}\frac{\left(  m+n-j\right)  !}{\left(  m-j\right)  !\left(  n-j\right)
!j!}\left(  -1\right)  ^{j}\left(  z\overline{z}\right)  ^{n-j}\quad (m<n),
\]
and satisfy the orthogonality relations%
\[
\frac{1}{\pi}\int\int_{\left\vert z\right\vert <1}Z_{mn}\left(  z,\overline
{z}\right)  \overline{Z_{kl}\left(  z,\overline{z}\right)  }dm_{2}\left(
z\right)  =\frac{\delta_{mk}\delta_{nl}}{m+n+1}.
\]
Suppose $f\left(  z,\overline{z}\right)  $ is continuous on the disk and has
coefficients%
\[
\widehat{f}_{mn}:=\int\int_{\left\vert z\right\vert <1}f\left(  z,\overline
{z}\right)  \overline{Z_{mn}\left(  z,\overline{z}\right)  }dm_{2}\left(
z\right)  ,m,n=0,1,2,\ldots
\]
then%
\[
f=\frac{1}{\pi}\sum_{m=0}^{\infty}\sum_{n=0}^{\infty}\left(  m+n+1\right)
\widehat{f}_{mn}Z_{mn},
\]
with convergence at least in the $L^{2}$-sense. As is typical of Fourier
expansions, the convergence behaviour is better for smoother functions $f$. If
$f$ is real then $\widehat{f}_{mn}=\overline{\widehat{f}_{nm}}$.
{\\[.5cm]}
Suppose $A$ is an $N\times N$ matrix whose numerical range is contained in the
unit disk (otherwise work with $B=c_{1}A+c_{0}I$, with $c_{1}>0$ so that
$\trace(B)=0$ and the range of $B$ satisfies the boundedness condition). We
may use Proposition 5.3 to
determine the moments of the shadow $P_{A}$ (and we write
$dP_{A}\left(  z\right)  =p_{A}\left(  z\right)  dm_{2}\left(  z\right)  $, so
that $p_{A}$ is the density). Thus
\[
\int_{\left\vert z\right\vert \leq1}z^{m}\overline{z}^{n}dP_{A}\left(
z\right)  =\frac{m!n!}{\left(  N\right)  _{m+n}}\left[  s^{m}t^{n}\right]
\xi_{A}\left(  s,t\right)  ^{-1},
\]
where $\xi_{A}\left(  s,t\right)  =\det\left(  I-sA-tA^*\right)  $ and
$\left[  s^{m}t^{n}\right]  g\left(  s,t\right)  $ denotes the coefficient of
$s^{m}t^{n}$ in the power series expansion of $g$ centered at $\left(
s,t\right)  =\left(  0,0\right)  $.
{\\[.5cm]}
It is then straightforward to compute the
Zernike coefficients of the density:%
\begin{eqnarray*}
\left(  p_{A}\right)  _{mn}^{\symbol{94}} 
&=& \int_{\left\vert z\right\vert \leq1}\overline{Z_{mn}
  \left(  z,\overline{z}\right)  }p_{A}\left(  z\right) dm_{2}\left(  z\right)
\\
&=& \int_{\left\vert z\right\vert \leq1}\sum_{j=0}^{n}
  \frac{\left(m+n-j\right)  !}{\left(  m-j\right)  !\left(  n-j\right)  !j!}
  \left(-1\right)  ^{j}z^{n-j}\overline{z}^{m-j}dP_{A}\left(  z\right)
\\
&=& \sum_{j=0}^{n}\frac{\left(  m+n-j\right)  !}{j!\left(  N\right)  _{m+n-2j}}
  \left(  -1\right)  ^{j}\left[  s^{n-j}t^{m-j}\right]  \xi_{A}
  \left(s,t\right)^{-1},
\end{eqnarray*}
for $m\geq n$ (and $\left(  p_{A}\right)  _{nm}^{\symbol{94}}=\overline
{\left(  p_{A}\right)  _{mn}^{\symbol{94}}}$ ). As an approximation, one may compute
$\left(  p_{A}\right)  _{mn}^{\symbol{94}}$ for all $\left(  m,n\right)  $
with $m+n\leq M$ for some $M$ (say 10 or 20). To write our formulas in
real terms with $z=x+iy$ (and $x^{2}+y^{2}\leq1$), let
\[
Q_{mn}\left(  u\right)  =\sum_{j=0}^{\min\left(  m,n\right)  }\frac{\left(
m+n-j\right)  !}{\left(  m-j\right)  !\left(  n-j\right)  !j!}\left(
-1\right)  ^{j}u^{\min\left(  m,n\right)  -j}.
\]
Note the trivial identity $ab+\overline{a}\overline{b}=2\mbox{Re}%
a\mbox{Re}b-2\mbox{Im}a\mbox{Im}b$. Thus the partial
sum for $p_{A}$ can be written as:%
\[
 \sum_{j=0}^{\left\lfloor M/2\right\rfloor }\left(  2j+1\right)  \left(
p_{A}\right)  _{jj}^{\symbol{94}}Q_{jj}\left(  x^{2}+y^{2}\right)  +
\]\[
 +2\sum_{j=1}^{M}\left(  \mbox{Re}\left(  \left(  x+\mathrm{i}%
y\right)  ^{j}\right)  \sum_{n=0}^{\left\lfloor \left(  M-j\right)
/2\right\rfloor }\left(  2n+j+1\right)  \mbox{Re}\left(  \left(
p_{A}\right)  _{n+j,n}^{\symbol{94}}\right)  Q_{n+j,n}\left(  x^{2}%
+y^{2}\right)  \right)
\]\[
 -2\sum_{j=1}^{M}\left(  \mbox{Im}\left(  \left(  x+\mathrm{i}%
y\right)  ^{j}\right)  \sum_{n=0}^{\left\lfloor \left(  M-j\right)
/2\right\rfloor }\left(  2n+j+1\right)  \mbox{Im}\left(  \left(
p_{A}\right)  _{n+j,n}^{\symbol{94}}\right)  Q_{n+j,n}\left(  x^{2}%
+y^{2}\right)  \right)  .
\]
(The factor $\pi$ has been ignored; it is merely a change of scale). It is a
matter for experimentation to produce useful graphs for a given matrix. The
polynomials tend to wiggle close to the edge of the disk; the graphs can not
be expected to precisely show the boundary of the numerical range, but they do
indicate the behaviour of $p_{A}$ in the interior.

\section{Numerical shadows and the higher--rank numerical ranges}

The rank--$k$
numerical ranges, denoted below by $\Lambda_k$, were introduced c. 2006 by Choi, Kribs, and \.Zyczkowski
as a tool to handle compression problems in quantum information theory. Since then
their theory and applications have been advanced with remarkable enthusiasm. The
sequence of papers [CHK\.Z2007,CGHK2008,W2008,LS2008], for example, led to a striking
extension of the classical Toeplitz--Hausdorff theorem (convexity of $W(M)$): \textbf{all} the
$\Lambda_k(M)$ are convex (though some may be empty), and they are intersections
of conveniently computable half--planes in $\dC$.
Among the many more recent papers concerning the $\Lambda_k(M)$, let us mention [LPS2009,GLW2010].
{\\[.5cm]}
Given a matrix $M\in M_N$ and $k\geq1$, Choi, Kribs, and \.Zyczkowski (see [CK\.Z2006a,CK\.Z2006b])
defined the rank--$k$ numerical range of $M$ as
\[
\Lambda_k(M)=\{\ld\in\dC:\exists P\in P_k\mbox{ such that }PMP=\ld P\},
\]
where $P_k$ denotes the set of rank--$k$ orthogonal projections in $M_N$. It is not hard
to verify that $\Lambda_K(M)$ can also be described as the set of complex $\ld$ such that
there is some $k$--dimensional  subspace $S$ of $\dC^N$ such that $(Mu,u)=\ld$ for \textbf{all}
unit vectors in $S$. In particular, we see that
\[
W(M)=\Lambda_1(M)\supseteq\Lambda_2(M)\supseteq\Lambda_3(M)\supseteq\dots\quad.
\]

\begin{figure}[ht!]
\centering
\subfloat[][
Shadow of $\diag(0,1,3,5)$.
] {%
\includegraphics[width=\figurewidth]{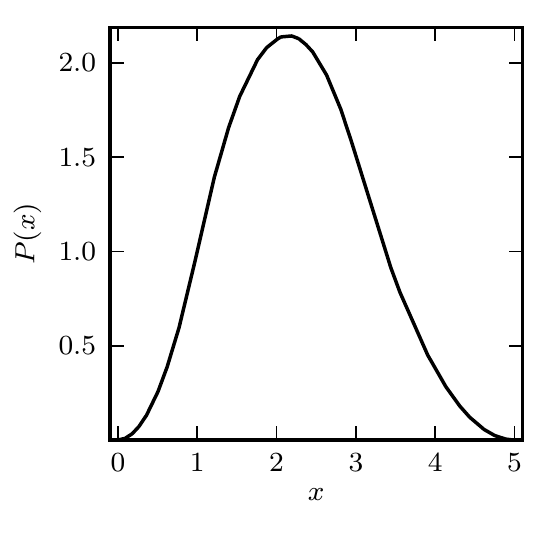}%
} 
\subfloat[][
Shadow of $\diag(0,1,3,5)$ with respect to real vectors.
] {%
\includegraphics[width=\figurewidth]{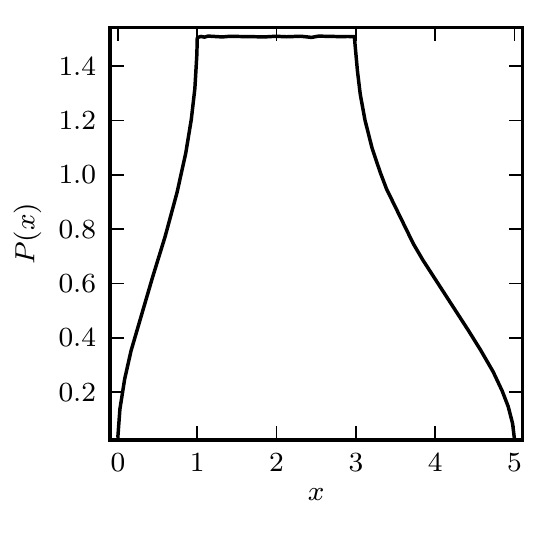}%
}%
\caption{ }
\label{fig1sec12}
\end{figure}
\begin{figure}[hb!]
\centering
\subfloat[][
Shadow of unitary matrix in $M_5$.
] {%
\includegraphics[width=\figurewidth]{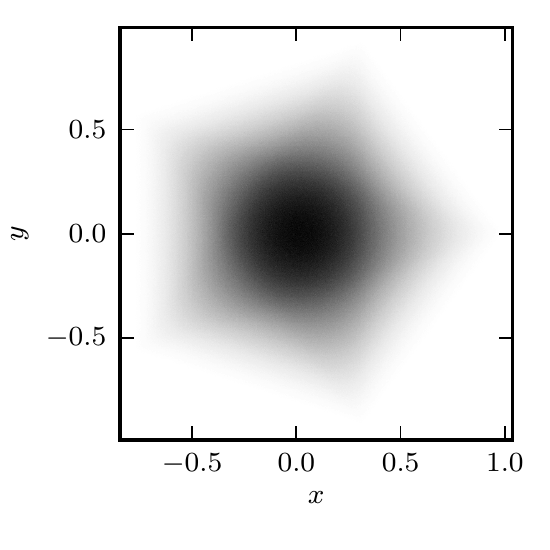}%
}%
\subfloat[][
Shadow of unitary matrix in $M_5$ with respect to real vectors.
] {%
\includegraphics[width=\figurewidth]{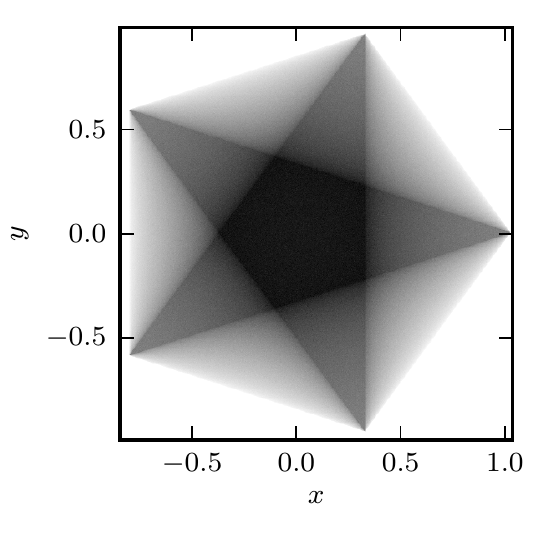}%
}%
\caption{ }
\label{fig2sec12}
\end{figure}
This point of view also suggests that these higher--rank numerical ranges should be visible as regions
of higher density within the numerical shadow of $M$. This idea is borne out, to some degree, by examining
shadow densities of various matrices. In Figure \ref{fig1sec12}(a), for example, we see the (one--dimensional) shadow
density of the Hermitian $\diag(0,1,3,5)$ - a spline of degree 2.  Here it is known that $\Lambda_2(M)=[1,3]$;
while the density is unimodal, there are values in [3,5] that are greater than some in [1,3]. If \textbf{real}
unit vectors are used in such experiments, the higher--rank numerical ranges often seem to be revealed more
clearly.; compare Figure \ref{fig1sec12}(b).


A similar phenomenon is seen in Figure \ref{fig2sec12}. Here $M$ is a unitary matrix in $M_5$ and it is known that
$\Lambda_2(M)$ is the inner pentagon (with interior) formed by lines joining the non--adjacent eigenvalues
in pairs. The shadow density is unimodal, but $\Lambda_2(M)$ is only seen distinctly in Figure \ref{fig2sec12}(b), where
only real unit vectors in $u\in\Omega_5$ are use to generate the values $(Mu,u)$.
The distinction between shadows based on complex vs real unit vectors is a reflection of the fact that the latter follow
a Dirichlet distribution (with parameter 1/2) rather than our usual measure $\mu$. See [B\.Z2006].

\ifthenelse{\equal{\LAA}{true}}
{
\section{References}
}
{ 
}




\begin{thebibliography}{99}

\bibitem[A1971]{A1971} T. W. Anderson, The Statistical Analysis of Time Series,
   Wiley, New York, 1971.

\bibitem[B1997]{B1997} R. Bhatia, Matrix Analysis, Springer--Verlag, New York, 1997.

\bibitem[B\.Z2006]{BZ2006} I. Bengtsson and K. \.Zyczkowski, Geometry of Quantum States, Cambridge UP, Cambridge, 2006. 

\bibitem[CK1985]{CK1985} E. W. Cheney and D. R. Kincaid, Numerical Mathematics and Computing, Brooks/Cole 1985.

\bibitem[CK\.Z2006a]{CKZ2006a} M.--D. Choi, D. Kribs, and K. \.Zyczkowski, Quantum error correcting codes from 
the compression formalism, Rep. Math. Phys. 58 (2006) 77--86.

\bibitem[CK\.Z2006b]{CKZ2006b} M.--D. Choi, D. Kribs, and K. \.Zyczkowski, Higher--rank numerical ranges and compression
problems, Linear Alg. Appl. 418 (2006) 828--839.

\bibitem[CHK\.Z2007]{CHKZ2007} M.--D. Choi, J. Holbrook, D. Kribs, and K. \.Zyczkowski, Higher--rank numerical ranges of unitary and normal matrices,
 Operators and Matrices 1 (2007) 409--426.

\bibitem[CGHK2008]{CGHK2008} M.--D. Choi, M. Giesinger, J. Holbrook, and D. Kribs, Geometry of higher--rank numerical ranges, 
Linear and Multilinear Algebra 56 (2008)  53-64.

\bibitem[DH1988]{DH1988} K. R. Davidson and J. Holbrook, Numerical radii of zero-one matrices, Michigan Math. J. 35 (1988) 261--267.

\bibitem[DJ2007]{DJ2007} D. \v Z. Djokovi\'c and C. R. Johnson, Unitarily achievable zero patterns and traces of words in $A$ and $A^*$,
Linear Alg. Appl. 421 (2007) 63--68.

\bibitem[D1971]{D1971} C. Davis, The Toeplitz--Hausdorff theorem explained, 
Canad. Math. Bull. 14 (1971) 245--246.

\bibitem[dB1976]{dB1976} C. de Boor,
Splines as linear combinations of B-splines, pp. 1-47
in Approximation Theory II (G.G. Lorentz, C. K. Chui, and L. L. Schumaker, eds.),
 Academic Press, New York, 1976.

\bibitem[GLW2010]{GLW2010} H.--L. Gau, C.--K. Li, and P. Y. Wu, 
Higher --rank numerical ranges and dilations,
J. Operator Theory 63 (2010) 181-189.

\bibitem[GR1997]{GR1997} K. E. Gustafson and D. K. M. Rao,
 Numerical Range, Springer, 1997.

\bibitem[GS2010]{GS2010} T. Gallay and D. Serre, The numerical measure of a complex matrix, arXiv:1009.1522v1 [math.FA], 8Sep2010

\bibitem[He1984]{He1984} S. Helgason, 
Groups and Geometric Analysis, Academic Press, New York, 1984.

\bibitem[Hn1974]{Hn1974} P. Henrici,
Applied and Complex Analysis, Vol. 1,
  Wiley-Interscience, New York 1974.

\bibitem[H2010]{H2010} J. Holbrook, Diagonal compressions of matrices and numerical shadows, colloquium, Feb 19, U. of Hawaii, 2010

\bibitem[HS2010]{HS2010} J. Holbrook and J.--P. Schoch, Theory vs. experiment: multiplicative inequalities for the numerical radius of commuting matrices, 
Operator Theory: Advances and Applications 202 (2010) 273--284.

\bibitem[JAG1998]{JAG1998} E. Jonckheere, F. Ahmad, and E. Gutkin, 
Differential topology of numerical range,
 Linear Alg. Appl. 279 (1998) 227--254.

\bibitem[L2001]{L2001} C.--K. Li, A survey on linear preservers of numerical ranges and radii, Taiwanese J. Math. 5 (2001) 477--496.

\bibitem[LPS2009]{LPS2009} C.--K. Li, Y.--T. Poon, and N.--S. Sze, Condition for the higher--rank numerical range to be non--empty, 
Linear and Multilinear Algebra 57 (2009) 365--368.

\bibitem[LS2008]{LS2008} C.--K. Li and N.--S. Sze, Canonical forms, higher rank numerical ranges, totally isotropic subspaces, and matrix equations, 
Proc. Amer. Math. Soc. 136 (2008) 3013--3023.

\bibitem[Mac1995]{Mac1995} I. G. Macdonald, 
 Symmetric Functions and Hall Polynomials,
 II ed., Clarendon Press, Oxford, 1995.

\bibitem[N1982]{N1982} K.--C. Ng, Some properties of doubly-stochastic matrices and distribution of density on a numerical range,
MPhil thesis, U of Hong Kong, 1982

\bibitem[S1940]{S1940} W. Specht, Zur Theorie der Matrixen II, 
Jahresber. Deutsche Math. 50 (1940) 19--23.

\bibitem[W2008]{W2008} H. Woerdeman, The higher rank numerical range is convex, Linear and Multilinear Algebra 56 (2008) 65--67.

 \bibitem[\.Z1999]{Z1999} K. {\.Z}yczkowski,
 Volume of the set of separable states II,
 Phys. Rev. A60 (1999) 3496--3507.

\bibitem[\.Z2009]{Z2009} K. \.Zyczkowski (with M.--D. Choi, C. Dunkl, J. Holbrook, P. Gawron,
J.Miszczak, Z. Puchala, and L. Skowronek), Generalized numerical range as a versatile tool to study quantum entanglement,
Oberwolfach, Dec 2009 (see Oberwolfach Report No. 59/2009, 34--37), 2009

\end{thebibliography}
\end{document}